\documentclass{amsart}

\title{On $12$-Congruences of Elliptic Curves}

\author{Sam Frengley} 
\address{University of Cambridge, DPMMS, Centre for Mathematical Sciences, Wilberforce Road, Cambridge CB3 0WB, UK}
\email{sam.frengley@gmail.com}

\date{September 14, 2023}

\usepackage{fullpage}		
\usepackage{colonequals}	
\usepackage{amsmath}		
\usepackage{amssymb}		
\usepackage{amsthm}			
\usepackage{amsfonts}		
\usepackage{enumitem}		
\usepackage[pagebackref=true,hypertexnames=false]{hyperref}

\usepackage{mathrsfs}
\usepackage{graphicx}
\usepackage{caption}
\usepackage{mathtools}
\usepackage{booktabs}
\usepackage{listings}
\usepackage{color}
\usepackage{float}
\usepackage{etoolbox}
\usepackage[utf8]{inputenc}

\usepackage{theoremref}
\newcommand{\threfpart}[2]{\thref{#1}\ref{#1_#2}}

\usepackage{tikz}
\usepackage{tikz-cd}

\usepackage{array}
\usepackage{makecell}
\usepackage{adjustbox}

\newcommand{\bbP}{\mathbb{P}}
\newcommand{\bbQ}{\mathbb{Q}}
\newcommand{\bbZ}{\mathbb{Z}}
\newcommand{\bbA}{\mathbb{A}}
\newcommand{\bbF}{\mathbb{F}}

\newcommand{\overbar}[1]{\mkern 1.5mu\overline{\mkern-1.5mu#1\mkern-1.5mu}\mkern 1.5mu}

\newcommand{\Kbar}{\overbar{K}}
\newcommand{\Qbar}{\overbar{\mathbb{Q}}}

\newcommand{\LMFDBLabel}[1]{\href{https://www.lmfdb.org/EllipticCurve/Q/#1/}{\texttt{#1}}}

\DeclareMathOperator{\Image}{im}

\DeclareMathOperator{\Gal}{Gal}
\DeclareMathOperator{\Aut}{Aut}
\DeclareMathOperator{\End}{End}

\DeclareMathOperator{\GL}{GL}
\DeclareMathOperator{\SL}{SL}

\DeclareMathOperator{\PGL}{PGL}
\DeclareMathOperator{\Proj}{Proj}
\DeclareMathOperator{\Spec}{Spec}
\DeclareMathOperator{\Bir}{Bir}

\newcommand{\eniii}{\textnormal{(\roman*)}}

\newtheorem{theorem}{Theorem}
\numberwithin{theorem}{section}

\newtheorem{lemma}[theorem]{Lemma}

\newtheorem*{theorem*}{Theorem}
\newtheorem*{lemma*}{Lemma}
\newtheorem*{prop*}{Proposition}
\newtheorem*{coro*}{Corollary}  
\theoremstyle{definition}\newtheorem{defn}[theorem]{Definition}
\theoremstyle{definition}\newtheorem{remark}[theorem]{Remark}
\theoremstyle{definition}\newtheorem{example}[theorem]{Example}
\theoremstyle{definition}\newtheorem*{defn*}{Definition}

\begin{document}

\begin{abstract}
	We construct infinite families of pairs of (geometrically non-isogenous) elliptic curves defined over $\bbQ$ with $12$-torsion subgroups that are isomorphic as Galois modules. This extends previous work of Chen~\cite{C_COEC} and Fisher~\cite{F_EMSFCOEC} where it is assumed that the underlying isomorphism of $12$-torsion subgroups respects the Weil pairing. Our approach is to compute explicit birational models for the modular diagonal quotient surfaces which parametrise such pairs of elliptic curves.
	
	A key ingredient in the proof is to construct simple (algebraic) conditions for the $2$, $3$, or $4$-torsion subgroups of a pair of elliptic curves to be isomorphic as Galois modules. These conditions are given in terms of the $j$-invariants of the pair of elliptic curves.
\end{abstract}

\maketitle

\vspace{0.2em}

\section{Introduction}
	Let $N \geq 1$ be an integer, $K$ be a perfect field of characteristic coprime to $N$, and $G_K \colonequals \Gal(\Kbar/K)$ be the absolute Galois group of $K$. Let $E/K$ and $E'/K$ be elliptic curves and let $E[N]$ and $E'[N]$ be their $N$-torsion subgroups. 
	
	We say that $E$ and $E'$ are \emph{$N$-congruent} (or \emph{$N$-congruent over $K$}) if there exists a $G_K$-equivariant isomorphism $\phi \colon E[N] \to E'[N]$. In this case the isomorphism $\phi$ is said to be an \emph{$N$-congruence}.

	Frey and Mazur conjectured that there exists a constant $C_0 > 1$ such that for any $N \geq C_0$ there are no pairs of non-isogenous $N$-congruent elliptic curves over $\bbQ$. It is necessary to exclude isogenies since for each integer $m$ coprime to $N$ the restriction of an $m$-isogeny to $E[N]$ is an $N$-congruence.

	Note that if $\phi \colon E[N] \to E'[N]$ is an $N$-congruence then there exists an element $r \in (\bbZ/N\bbZ)^\times$ such that for all $P, Q \in E[N]$ we have
	\begin{equation*}
		e_{E, N}(P, Q) = e_{E', N}(\phi(P), \phi(Q))^r
	\end{equation*}
	where $e_{E, N}$ and $e_{E', N}$ denote the Weil pairings on $E$ and $E'$ respectively. In this case we say that $\phi$ is an \emph{$(N, r)$-congruence} or an \emph{$N$-congruence with power $r$}. For each integer $m$ coprime to $N$ composing $\phi$ with multiplication by $m$ changes the power of $\phi$ by a factor of $m^2$, so we consider the power of $\phi$ only up to multiplication by a square in $(\bbZ/N\bbZ)^\times$.
	
	We are concerned with the following question: for which pairs $(N, r)$ do there exist examples (or infinite families of examples) of pairs of non-isogenous $(N, r)$-congruent elliptic curves over $\bbQ$?
	
	In the case when $N = p$ is a prime number, infinite families of geometrically non-isogenous $(p,r)$-congruent elliptic curves over $\bbQ$ have been found for every power $r$ when $p \leq 13$ (see Table~\ref{table:Nr-congs-refs}). A pair of non-isogenous $(17, 3)$-congruent elliptic curves appearing in the LMFDB~\cite{lmfdb} has been found by Cremona (see \cite{B_OSRCBTEC}, \cite{CF_GMFTSTOCBEC}, or \cite{F_OFO7A11CEC}). Fisher~\cite{F_OPO17CEC} also found a pair of non-isogenous $(17, 1)$-congruent elliptic curves over $\bbQ$, and has conjectured that these are the only $17$-congruences over $\bbQ$. Moreover Fisher conjectured that there exist no pairs of non-isogenous, $p$-congruent, elliptic curves over $\bbQ$ for $p \geq 19$ (see \cite[Conjecture~1.1]{F_OPO17CEC}). 

	For composite $N$ infinite families of geometrically non-isogenous $(N, r)$-congruent elliptic curves have been found for every power $r$ when $N \leq 10$ and for $(N, r) = (12, 1)$ (see Table~\ref{table:Nr-congs-refs}).
	
	On the other hand in \cite{F_COECAFNSMNGR} we found examples of elliptic curves over $\bbQ$ admitting an $N$-congruence with a (non-isogenous) quadratic twist for each even $N \leq 24$, $N = 28$, $30$, $32$, $36$, and $48$. These examples arise in infinite families except in the cases when $N = 30$, $32$, or $48$. When $N = 12$ we showed that there are infinite families of elliptic curves admitting $(12, 7)$ and $(12, 11)$-congruences with a (non-isogenous) quadratic twist. However, in these examples the elliptic curves in question are geometrically isogenous and their mod $12$ Galois representations are not surjective.

  \begin{table}[t]
    \centering
    \begin{tabular}{c|>{\centering\arraybackslash}p{2.5cm}||c|>{\centering\arraybackslash}p{2.5cm}||c|>{\centering\arraybackslash}p{2.5cm}}
      $(N,r)$ & References                                                           & $(N,r)$ & References                                                                & $(N,r)$  & References                                   \\
      \hline
      $(2,1)$ & \cite{F_THOAG1C,K_HMSFSDAESOG2FF,RS_M2ROEC}                          & $(6,5)$ & \cite{C_COEC,K_HMSFSDAESOG2FF}                                            & $(9,2)$  & \cite{F_OFO9CEC,F_EMSFCOEC,K_HMSFSDAESOG2FF} \\
      $(3,1)$ & \cite{ATT_HPA3TS,F_THOAG1C,K_CFOECWPM3RVHACC,RS_FOECWCMPR}           & $(7,1)$ & \cite{C_COEC,F_OFO7A11CEC,HK_SLCMXE7}                                     & $(10,1)$ & \cite{F_EMSFCOEC,K_HMSFSDAESOG2FF}           \\
      $(3,2)$ & \cite{F_THOAG1C,BHLS_G2CAJWAGNOP,K_HMSFSDAESOG2FF,K_CFOECWPM3RVHACC} & $(7,3)$ & \cite{C_COEC,F_OFO7A11CEC,F_EMSFCOEC,PSS_TOX7APSTX2Y3Z7,K_HMSFSDAESOG2FF} & $(10,3)$ & \cite{F_EMSFCOEC}                            \\
      $(4,1)$ & \cite{F_THOAG1C,S_EFOECWPMNR}                                        & $(8,1)$ & \cite{C_COEC,C_FOECWTSM8R}                                                & $(11,1)$ & \cite{F_OFO7A11CEC,F_EMSFCOEC,KR_MQOM11ROEC} \\
      $(4,3)$ & \cite{BD_TAOG2CW44SJ,F_THOAG1C,K_HMSFSDAESOG2FF}                     & $(8,3)$ & \cite{C_COEC,C_FOECWTSM8R,F_EMSFCOEC}                                     & $(11,2)$ & \cite{F_OFO7A11CEC,K_HMSFSDAESOG2FF}         \\
      $(5,1)$ & \cite{F_THOAG1C,K_HMSFSDAESOG2FF,RS_FOECWCMPR}                       & $(8,5)$ & \cite{C_COEC,C_FOECWTSM8R,F_EMSFCOEC}                                     & $(12,1)$ & \cite{C_COEC,F_EMSFCOEC}                     \\
      $(5,2)$ & \cite{F_THOAG1C}                                                     & $(8,7)$ & \cite{C_COEC,C_FOECWTSM8R,F_EMSFCOEC,K_HMSFSDAESOG2FF}                    & $(13,1)$ & \cite{F_OFO13CEC}                            \\
      $(6,1)$ & \cite{C_COEC,P_CEAM6TQUCED,R_EFOECWPM6R,RS_M6ROEC}                   & $(9,1)$ & \cite{F_OFO9CEC,F_EMSFCOEC}                                               & $(13,2)$ & \cite{F_OFO13CEC}                            \\
    \end{tabular}
    \caption{Pairs $(N,r)$ for which it is known that there are infinite families of geometrically non-isogenous $(N,r)$-congruent elliptic curves over $\bbQ$, together with authors who have computed a model for either of the moduli spaces $X_E^r(N)$ or $Z(N,r)$ defined in Section~\ref{sec:modularcurve}.} \label{table:Nr-congs-refs}
  \end{table}
  
	For any $N \geq 2$ and $r \in (\bbZ/N\bbZ)^\times$ let $Z(N, r)$ be the (coarse) moduli space which parametrises triples $(E, E', \phi)$ where $E$ and $E'$ are elliptic curves (defined up to simultaneous quadratic twist) and $\phi \colon E[N] \to E'[N]$ is an $(N, r)$-congruence (defined up to composition with an automorphism of $E$ or $E'$). Following Kani--Schanz~\cite{KS_MDQS} we call $Z(N, r)$ a  \emph{modular diagonal quotient surface}. The surface $Z(N, r)$ is naturally equipped with an involution which reverses the roles of $E$ and $E'$. We write $W(N, r)$ for the quotient of $Z(N, r)$ by this involution.

	The main contribution of this article is to give explicit equations for the surfaces $Z(12, r)$ for each $r \in (\bbZ/12\bbZ)^\times$. We use these equations to find infinite families of $(12, r)$-congruent elliptic curves which are not geometrically isogenous. 

	The surface $Z(12, 1)$ is an elliptic K3 surface and models were computed by Chen and Fisher (see \cite[Chapter~7]{C_COEC} and \cite{F_EMSFCOEC}) using equations for the twist $X_E^1(12)$ of $X(12)$ which parametrises elliptic curves which are $(12, 1)$-congruent to $E$ (see Section~\ref{sec:twistX_E}). Chen~\cite[Chapter~7.4]{C_COEC} also computed equations for the twist $X_E^7(12)$ and it may be possible to compute a simple birational model for $Z(12, 7)$ using these. Using different techniques we prove:

	\begin{theorem} \thlabel{thm:Zequations}
		Each of the surfaces $W(12,r)$ are rational and the surfaces $Z(12, r)$ are birational over $\bbQ$ to the affine surfaces 
		\begin{equation*}
			z^2 = F_{12,r}(u,v)
		\end{equation*}
		where $F_{12,r}(u,v) \in \bbQ[u,v]$ are the polynomials
		\begin{flalign*}
			F_{12, 1}(u,v) =& -4 u^4 - 2 (2 v^2 + 3) u^3 - (v^4 - 13 v^2 + 3)u^2 + 2(5 v^4 + 9 v^2) u + v^6 - 6 v^4 + 9 v^2, &
		\end{flalign*}
		\begin{flalign*}
			F_{12, 5}(u,v) =& v^2 u^8 + 4 v^2 u^7 + (-4 v^4 - 20 v^2 + 4) u^6 + (-12 v^4 - 32 v^2 - 8) u^5 + 
			(6 v^6 + 44 v^4 + 160 v^2 - 76) u^4& \\&+ (12 v^6 + 64 v^4 - 64 v^2 + 352) u^3 +
			(-4 v^8 - 28 v^6 - 204 v^4 - 232 v^2 - 592) u^2 \\&+ (-4 v^8 - 32 v^6 + 72 v^4 
			+ 224 v^2 + 448) u + v^{10} + 4 v^8 + 40 v^6 + 52 v^4 - 32 v^2 - 128 	,	
		\end{flalign*}
		\begin{flalign*}
			F_{12, 7}(u,v) =& -(u - 1)((v^2 + 1) u + v^2 - 1)(u^4 + 6 u^3 + (-2 v^2 + 12) u^2 + (-6 v^2 + 8) u
			+ v^4 - v^2),		&
		\end{flalign*}
		\begin{flalign*}
			F_{12, 11}(u,v) =& -(u + 1)((v^4 + 6 v^2 + 1) u - 7 v^4 - 2 v^2 + 1)(v^4 u^4 + 8 v^4 u^3 + (-9 v^6 
			+ 24 v^4 + 11 v^2) u^2 &\\&+ (-54 v^6 - 36 v^4 - 6 v^2) u + 27 v^8 + 27 v^6 + 9 v^4 
			- v^2 - 1).
		\end{flalign*} 
		Moreover the double covers $Z(12, r) \to W(12, r)$ are given by $(u,v,z) \mapsto (u,v)$. 
	\end{theorem}

	The maps $Z(12, r) \to X(1) \times X(1)$ which give the moduli interpretations for the surfaces $Z(12, r)$ are too complicated to reproduce here but may be recovered from the computations in Section~\ref{sec:compute}. We record them in the electronic data corresponding to this article~\cite{ME_ELECTRONIC_Z12}. 
	
	\begin{remark}
		Note that the polynomials $F_{12, r}(u,v)$ contain only even powers of $v$.	In Section~\ref{sec:fibprod} we will show that $Z(12, r)$ is a double cover of a surface $Z^+(12,r)$. This double cover is the quotient by the involution $v \mapsto -v$.  
	\end{remark}

	\begin{remark}
		The model for the surface $Z(12, 1)$ given in \thref{thm:Zequations} admits a genus $1$ fibration over the $t$-line given by $y^2 = F_{12, 1}(x, t)$ (in fact, this fibration is elliptic, see \thref{ex:121}). 
		
		Since the surfaces $Z(12, r)$ are of general type when $r \neq 1$ they cannot admit genus $1$ fibrations. However, they do admit genus $2$ fibrations over the $t$-line which are given by $y^2 = F_{12, 5}\left(\frac{-x + 4}{2t}, \frac{x}{2t}\right)$, $y^2 = F_{12, 7}(x, t)$, and $y^2 = F_{12, 11}(x, t)$ respectively.
	\end{remark}

	It is possible to find pairs of $12$-congruent elliptic curves by sieving for examples where both $E$ and $E'$ are contained in the LMFDB \cite{lmfdb} (as described in \cite[Section~3]{CF_GMFTSTOCBEC}). The power of these congruences may be determined using the equations of Rubin, Silverberg, and Fisher \cite{F_THOAG1C,RS_FOECWCMPR,RS_M2ROEC,S_EFOECWPMNR} for $X_E^r(N)$ when $N = 3$ and $4$. For example the pairs of elliptic curves with LMFDB labels $(\LMFDBLabel{55.a1}, \LMFDBLabel{1045.b1})$, $(\LMFDBLabel{60450.cx2}, \LMFDBLabel{60450.cw2})$, and $(\LMFDBLabel{735.d2}, \LMFDBLabel{24255.bh1})$ are not geometrically isogenous and are $(12, 1)$, $(12, 5)$, and $(12, 7)$-congruent respectively. 

	On the other hand we were not able to find any examples of $(12, 11)$-congruences between geometrically non-isogenous curves in the LMFDB. The rational point $\left( \frac{-19}{21}, \frac{3}{7}, \frac{1844480}{2470629} \right)$ on the model for $Z(12, 11)$ in \thref{thm:Zequations} gives rise to the pair of geometrically non-isogenous $(12, 11)$-congruent elliptic curves
	\begin{align*}
		E  &: y^2 + xy = x^3 - 21666120x -57035036608, \\
		E' &: y^2 + xy = x^3 + 398520965x + 166506419482597
	\end{align*}
	which have conductor $4\,976\, 690$. Moreover, among pairs of geometrically non-isogenous $(12, 11)$-congruent elliptic curves which we found, these curves have the smallest conductor (see Section~\ref{sec:smallconductor}).

	By exhibiting a rational curve on the elliptic K3 surface $Z(12,1)$ Chen found an infinite family of $(12,1)$-congruent elliptic curves over $\bbQ$ \cite[Proposition~1.7.11]{C_COEC}. Fisher gave a fibration with a section of infinite order, thus giving an infinite family of $(12, 1)$-congruent elliptic curves over $\bbQ(t)$ \cite[Corollary~1.3]{F_EMSFCOEC}.
	
	In contrast Kani--Schanz~\cite{KS_MDQS} showed that when $r \neq 1$ the surfaces $Z(12, r)$ are of general type. In this case we expect (by the Bombieri--Lang conjecture) that there are at most finitely many curves of geometric genus $0$ or $1$ on $Z(12, r)$. 

	By finding curves of genus $0$ or $1$ on the models for $Z(12, r)$ in \thref{thm:Zequations} we deduce:

	\begin{theorem} \thlabel{thm:infmany}
		For each $r \in (\bbZ/12\bbZ)^\times$ there are infinitely many pairs of $j$-invariants of $(12, r)$-congruent elliptic curves $E/\bbQ$ and $E'/\bbQ$ which are not geometrically isogenous. 
		
		Moreover there are infinite subfamilies such that the mod $12$ Galois representation attached to $E/\bbQ$ (and $E'/\bbQ$) is surjective.
	\end{theorem}

	\begin{remark}
		Since a pair of elliptic curves are $(12, r)$-congruent if and only if they are both $(3,r)$ and $(4,r)$-congruent, the families of elliptic curves in \thref{thm:infmany} may be checked using known formulae for the twists $X_E^r(N)$ when $N = 3, 4$. Such equations have been computed by Rubin and Silverberg \cite{RS_FOECWCMPR,RS_M2ROEC,S_EFOECWPMNR} when $r = 1$, and by Fisher~\cite{F_THOAG1C} for each $r$.
	\end{remark}

	\subsection{Outline of the paper} In Section~\ref{sec:modularcurve} we give background on the invariant theory of the modular curve $X(N)$ when $N \leq 4$ following Klein~\cite{K_LOTIASOEOTFD} and Fisher~\cite{F_THOAG1C}. We describe $Z(N, r)$ as a quotient of $X(N) \times X(N)$ following Kani and Schanz~\cite{KS_MDQS}.
	
	In Section~\ref{sec:234} we use invariant theoretic techniques to find simple (symmetric) conditions for a pair of elliptic curves to be $2$, $3$, or $4$-congruent (note that our condition for $2$-congruence was already proved by Fisher, see \cite[Lemma~3.5]{F_EMSFCOEC}).  
	
	We describe the relationship between the surface $Z(N_1 N_2, r)$ and the fibre product $Z(N_1, r) \times_{Z(1)} Z(N_2, r)$ when $N_1$ and $N_2$ are coprime integers in Section~\ref{sec:fibprod}. We employ this description to compute explicit models for the surfaces $Z(12, r)$ and prove \thref{thm:Zequations} in Section~\ref{sec:compute}. 
	
	Finally in Section~\ref{sec:curves} we find ``interesting'' curves on the surfaces $Z(12, r)$. We prove \thref{thm:infmany} in Section~\ref{sec:inffams} by finding curves with infinitely many rational points (and which do not arise as Hecke correspondences) on the models for $Z(12, r)$ in \thref{thm:Zequations}. More precisely we give explicit examples of pairs of geometrically non-isogenous $(12, r)$-congruent elliptic curves which are defined over $\bbQ(t)$ when $r = 1$, $7$, and $11$ (see Examples~\ref{ex:121}, \ref{ex:127}, and \ref{ex:1211}). When $r = 5$ we give a family defined over $\bbQ(\mathcal{C})$ where $\mathcal{C}/\bbQ$ is an elliptic curve of rank $1$ (see \thref{ex:125}). 
	
	\subsection{Acknowledgements} 
		I would like to thank my supervisor Tom Fisher for suggesting this topic and for many insightful conversations and comments on earlier versions of this article. I would also like to express my gratitude to Jef Laga and the anonymous referee for several comments on earlier versions of this paper and to Noah Porcelli, Tony Scholl, and Alice Silverberg for useful discussions and suggestions. I am grateful to the Woolf Fisher and Cambridge Trusts for their financial support.

\section{The modular curve \texorpdfstring{$X(N)$}{X(N)} and the surface \texorpdfstring{$Z(N, r)$}{Z(N, r)}} \label{sec:modularcurve}
	Let $N \geq 1$ be an integer, and let $K$ be a field of characteristic coprime to $6N$.

	Consider the $G_K$-module $\mu_N \times \bbZ/N\bbZ$ endowed with the structure of symplectic abelian group via the pairing $\langle (\zeta, c), (\xi, d) \rangle = \zeta^d \xi^{-c}$. Let $\Gamma\!_N$ denote the group of symplectic automorphisms of $\mu_N \times \bbZ/N\bbZ$, and note that $\Gamma\!_N$ is isomorphic to $\SL_2(\bbZ/N\bbZ)$ as an abstract group. 

	Let $Y(N)/K$ be the (non-compact) modular curve of full level $N$. For each field $L/K$, the $L$-rational points on $Y(N)$ parametrise pairs $(E/L, \iota)$ where $E/L$ is an elliptic curve and $\iota \colon E[N] \to \mu_N \times \bbZ/N\bbZ$ is a $G_L$-equivariant isomorphism which is symplectic with respect to the Weil pairing on $E[N]$. Recall that when $N = 1, 2$ the moduli space $Y(N)$ is coarse, and the pair $(E/L, \iota)$ is only determined up to twist.
	
	An element $g \in \Gamma\!_N$ acts on $Y(N)$ via $(E, \iota) \mapsto (E, g \iota)$ and this action extends uniquely to an automorphism of $X(N)$. Elements $g, g' \in \Gamma\!_N$ induce the same action on $X(N)$ if and only if they differ by an automorphism of $E$. In particular, this action descends to an action of $\Gamma\!_N/\{\pm 1\}$. 

	\subsection{The Twist \texorpdfstring{$X_E^r(N)$}{XE{\textasciicircum}r(N)}} \label{sec:twistX_E}
		Let $E/K$ be an elliptic curve and $r$ be an integer coprime to $N$. Consider a $\overbar{K}$-isomorphism $\iota_r \colon E[N] \to \mu_N \times \bbZ/N\bbZ$ such that $e_N(P, Q)^r =  \langle \iota_r(P), \iota_r(Q) \rangle$. By the twisting principle $\iota_r$ gives rise to a Galois cohomology class $[\iota_r] \in H^1(K, \Gamma\!_N)$. 

		The $G_K$-equivariant homomorphism $\Gamma\!_N \to \Aut(X(N)_{\overbar{K}})$ induces a map on Galois cohomology $H^1(K, \Gamma\!_N) \to H^1(K, \Aut(X(N)_{\overbar{K}}))$. Again by the twisting principle we may associate to the image of $[\iota_r]$ a twist $X_E^r(N)/K$ of $X(N)$. 
		
		For each field $L/K$, the $L$-rational points on the curve $X_E^r(N)$ parametrise pairs $(E'/L, \phi)$ where $\phi \colon E[N] \to E'[N]$ is an $(N, r)$-congruence. 

	\subsection{The Surface \texorpdfstring{$Z(N, r)$}{Z(N, r)}}
		Recall that for each $N \geq 2$ and $r \in (\bbZ/N\bbZ)^\times$ we let $Z(N, r)/K$ be the (coarse) moduli space of triples $(E, E', \phi)$ where $E/K$ and $E'/K$ are elliptic curves (defined up to simultaneous quadratic twist) and $\phi \colon E[N] \to E'[N]$ is an $(N, r)$-congruence. Further recall that we write $W(N, r)$ for the quotient of $Z(N, r)$ by the involution reversing the roles of $E$ and $E'$.
		
		We will write $Z(1)$ for the surface $Y(1) \times Y(1)$ and $W(1)$ for the quotient of $Z(1)$ by the natural involution swapping the roles of $j$ and $j'$. It will be useful for us to note that $W(1)$ is birational to $\bbA^2$ and that we may take the quotient $Z(1) \to W(1)$ to be given by $(j, j') \mapsto (jj', (j-1728)(j'-1728))$.
		
		We then have a commutative diagram:
		\[\begin{tikzcd}
			{Z(N, r)} & {W(N,r)}\\
			{Z(1)} & {W(1)}			
			\arrow[from=1-1, to=1-2]
			\arrow[from=1-1, to=2-1]
			\arrow[from=1-2, to=2-2]
			\arrow[from=2-1, to=2-2]
		\end{tikzcd}\]
		
		Define the group $\Delta_r$ to be the subgroup $\{ (g,  \varepsilon_r(g) ) : g \in \Gamma\!_N\} \subset \Gamma\!_N \times \Gamma\!_N$, where $\varepsilon_r$ denotes conjugation by the matrix $\big(\begin{smallmatrix} r & 0\\ 0 & 1 \end{smallmatrix}\big)$. 

		The following lemma is due to Kani--Schanz~\cite{KS_MDQS} and Fisher~\cite[Lemma~3.2]{F_OFO13CEC}.

		\begin{lemma}[{\cite[Lemma~3.2]{F_OFO13CEC}}, \cite{KS_MDQS}] \thlabel{lemma:MDQS}
			The surface $Z(N, r)$ is birational over $K$ to the quotient of $X(N) \times X(N)$ by $\Delta_r$.
		\end{lemma}
	
	\subsection{The Action of \texorpdfstring{$\Gamma\!_N$}{Gamma N} on \texorpdfstring{$X(N)$}{X(N)}} \label{sec:action}
		We briefly describe the invariant theory of the modular curves $X(N)$ when $N = 2$, $3$, and $4$. In these cases the invariant theory of $X(N)$ has been extensively studied, dating back to Klein~\cite{K_LOTIASOEOTFD}. We follow the more modern treatment of Fisher~\cite[Sections~3--4]{F_THOAG1C}.
		
		For $N = 2, 3, 4$ we identify the modular curve $X(N)$ with $\bbP^1$. We define the (finite) subgroups $G_N \subset \GL_2(\Kbar)$ to be the subgroups generated by
		\begin{align*}
			N = 2 &: \qquad 
							\begin{pmatrix}
								1/2 & 1/16 \\ 12 & -1/2
							\end{pmatrix}
							\quad\hspace{0.65em} \text{and} \quad\hspace{0.6em}	
							\begin{pmatrix}
								1 & 0 \\ 0 & -1
							\end{pmatrix}\\
			N = 3 &: \qquad 
							\frac{1}{\sqrt{-3}}\begin{pmatrix}
								1 & 1/3 \\ 6 & -1
							\end{pmatrix}
							\quad \text{and} \quad	
							\zeta_3 \begin{pmatrix}
								1 & 0 \\ 0 & \zeta_3
							\end{pmatrix}\\
			N = 4 &: \qquad 
							\frac{1}{\sqrt{-2}}\begin{pmatrix}
								1 & 1/2 \\ 2 & -1
							\end{pmatrix}
							\quad \text{and} \quad	
							\zeta_8^{-1} \begin{pmatrix}
								1 & 0 \\ 0 & \zeta_8^2
							\end{pmatrix}
		\end{align*}
		where $\zeta_N \in \Kbar$ is a primitive $N^{\text{th}}$ root of unity. 
		
		By mapping the generators $\big(\begin{smallmatrix} 0 & 1\\ -1 & 0 \end{smallmatrix}\big)$ and $\big(\begin{smallmatrix} 1 & 1\\ 0 & 1 \end{smallmatrix}\big)$ of $\Gamma\!_N$ to the above generators of (the projective image of) $G_N$ we obtain a $G_{K}$-equivariant homomorphism $\bar{\rho} \colon \Gamma\!_N \to \PGL_{2}(\Kbar)$ giving the action of $\Gamma\!_N$ on $X(N)$.

		Let $\widetilde{G}_N$ be the commutator subgroup of $G_N$. It is a classical theorem of Klein that for each $N = 2, 3, 4$ the invariant ring $K[x_0, x_1]^{\widetilde{G}_N}$ is generated by invariants $c_4$, $c_6$, and $D$ of weights $4N$, $6N$, and $12$ subject only to the relation $c_4^3 - c_6^2 = 1728 D^N$ (see \cite[Theorem 3.3]{F_THOAG1C}). Explicitly we take
		\begin{align*}
			N &= 2 : & 	D &= x_0(64 x_0^2 - x_1^2) 		& 	c_4 &= 192 x_0^2 + x_1^2  					& 	c_6 &= x_1 (576 x_0^2 - x_1^2) \\
			N &= 3 : & 	D &= -x_0(27x_0^3 + x_1^3) 		& 	c_4 &= -x_1(216 x_0^3 - x_1^3) 				& 	c_6 &= 5832 x_0^6 - 540 x_0^3 x_1^3 - x_1^6 \\
			N &= 4 : & 	D &= x_0 x_1 (16 x_0^4 - x_1^4) & 	c_4 &= 256 x_0^8 + 224 x_0^4 x_1^4 + x_1^8 	& 	c_6 &= 4096 x_0^{12} - 8448 x_0^8 x_1^4 - 528 x_0^4 x_1^8 + x_1^{12}.
		\end{align*}
		
		Moreover, by considering ramification points it can be shown that the map $X(N) \to X(1)$ given by $j = c_4^3/D^N$ is the $j$-map (see \cite[Lemma~4.4]{F_THOAG1C}).

		With the models chosen above the forgetful map $X(4) \to X(2)$ is given by $[x_0 : x_1] \mapsto [-x_0^2 x_1^2 : 16x_0^4 + x_1^4]$. 

\section{Symmetric equations for families of \texorpdfstring{$N$}{N}-congruent elliptic curves} \label{sec:234}

	Several approaches have been applied to parametrising (universal) families of $2$, $3$, and $4$-congruent elliptic curves. 
	
	Equations for $X_E^1(N)$ when $N = 2$, $3$, $4$, and $5$ were computed by Rubin and Silverberg (see \cite{RS_FOECWCMPR,RS_M2ROEC,S_EFOECWPMNR}). Using a direct approach Lario--Rio computed equations for $X_E^r(3)$ when $r = 1$ and $2$ (see \cite[Section~3]{LR_AOETEIB2K}).
	
	Fisher~\cite{F_THOAG1C} described how when $N = 2$, $3$, $4$, and $5$ the curve $X_E^r(N)$ is related to the Hessian of a genus one model of degree $N$, and hence computed explicit equations for $X_E^r(N)$ using classical invariant theory. When $N = 3$ similar approaches have been employed in \cite{ATT_HPA3TS} and \cite{K_CFOECWPM3RVHACC}. The former considers all powers, and the latter pays careful attention to the case when $r = 1$ and the characteristic of $K$ is $2$.
	
	Kumar~\cite{K_HMSFSDAESOG2FF} (for every $N \leq 11$) and Br{\"o}ker--Howe--Lauter--Stevenhagen~\cite{BHLS_G2CAJWAGNOP} (when $N = 3$), and Bruin--Doerksen~\cite{BD_TAOG2CW44SJ} (when $N = 4$) have computed explicit equations for the surface $Z(N, -1)$ using the description of $Z(N, -1)$ as the moduli space of genus $2$ curves with $(N,N)$-split Jacobians. The equations of Kumar and Br{\"o}ker--Howe--Lauter--Stevenhagen are also symmetrical in $E$ and $E'$.

	Our approach is to use \thref{lemma:MDQS} to compute models for the surfaces $Z(N, r)$ for each $(N, r) = (2, 1)$, $(3, 1)$, $(3, 2)$, $(4, 1)$, and $(4, 3)$. An advantage of our approach is that by viewing $E$ and $E'$ symmetrically we are able to exploit natural factorisations of the maps $Z(N, r) \to Z(1) = Y(1) \times Y(1)$ to provide simple conditions for $E$ and $E'$ to be $2$, $3$, or $4$-congruent (see \thref{thm:plys}). Indeed, it is this observation which makes it practical for us to compute (birational) models for $Z(12, r)$ via fibre products.

	\subsection{Computing quotients of \texorpdfstring{$X(N) \times X(N)$}{X(N) x X(N)}}
		We compute $Z(N,r)$ from its description as a quotient of $X(N) \times X(N)$ under the (twisted) diagonal action of $\Gamma\!_N$. The method we describe has been successfully employed by Fisher to compute birational models for $Z(13, r)$ and $Z(17, r)$ (see \cite{F_OFO13CEC} and \cite{F_OPO17CEC}). 

    For each $N = 2,3,4$ let $\Lambda_r \subset G_N \times G_N$ denote the preimage of $(\bar{\rho} \times \bar{\rho})(\Delta_r)$. Let $f(x_0, x_1, x_0', x_1')$ be a polynomial which is homogeneous of degrees $m$ and $m'$ in the two sets of variables.
		
		\begin{defn}
			We say that $f$ is an \emph{$r$--bi-invariant} if $f$ is invariant under the action of $\Lambda_r$ on the polynomial ring $\Kbar[x_0, x_1, x_0', x_1']$. If $f(x_0, x_1, x_0', x_1') = f(x_0', x_1', x_0, x_1)$
			we say that $f$ is \emph{symmetric}.
		\end{defn}

    The action of $\Delta_r$ on $X(N) \times X(N)$ induces an action on $Y(N) \times Y(N) = \Spec R$. By \thref{lemma:MDQS} the surface $Z(N,r)$ is birational to $\Spec R^{\Delta_r}$. To compute a birational model for $Z(N,r)$ it is enough to compute generators and relations for the $K$-algebra $R^{\Delta_r}$, together with their relations with the rational functions $j$ and $j'$, where $j$ (respectively $j'$) is the rational function giving the $j$-invariant on each copy of $Y(N)$. Similarly if $S \subset R^{\Delta_r}$ is the subalgebra of elements invariant under the natural involution on $Y(N) \times Y(N)$, then $W(N,r)$ is birational to $\Spec S$.  

    We do not attempt to fully describe the invariant ring $R^{\Delta_r}$. Instead we use the following observations:

    If $A \subset K[x_0,x_1,x_0',x_1']$ is the subalgebra of $r$--bi-invariants, then the inclusion $Y(N) \times Y(N) \subset X(N) \times X(N)$ induces a rational map $Z(N,r) \dashrightarrow \Spec R^{\Delta_r} \dashrightarrow \Proj A$. Moreover if $B \subset A$ is the subalgebra of symmetric $r$--bi-invariants, then we have a rational map $W(N,r) \dashrightarrow \Spec S \dashrightarrow \Proj B$.
	
		Consider non-zero $r$--bi-invariants $w_0, ..., w_k \in K[x_0,x_1,x_0',x_1']$ such that for each $i$ the $r$--bi-invariant $w_i$ is an element of the $K$-vector space of symmetric $r$--bi-invariants of homogeneous bi-degree $(m, m)$. Let $\varphi_{N,r}$ be the map $X(N) \times X(N) \dashrightarrow \bbP^k$ given by $([x_0 : x_1], [x_0' : x_1']) \mapsto [w_0 : ... : w_k]$. Since $w_0, ..., w_k$ are symmetric and invariant under the action of $\Lambda_r$ we obtain a factorisation  
		\[\begin{tikzcd}
			{X(N) \times X(N)} \\
			{W(N,r)} & {\Image(\varphi_{N,r}).}
			\arrow[from=1-1, to=2-1]
			\arrow[dashed, "\psi_{N,r}", from=2-1, to=2-2]
			\arrow[dashed, "\varphi_{N,r}", from=1-1, to=2-2]
		\end{tikzcd}\]
		
		When $N = 2$, $3$, and $4$ we will carefully choose symmetric $r$--bi-invariants $w_0, ..., w_k$ so that $\psi_{N,r}$ is birational.
		
	\subsection{Symmetric bi-invariants} \label{sec:biinvs}
		We now give choices of $r$--bi-invariants such that the map $\psi_{N,r}$ is birational for each $N = 2$, $3$, and $4$. We consider the action of $\Gamma\!_N$ on $X(N) \cong \bbP^1$ via the representation $\bar{\rho}$ in Section~\ref{sec:action}. 
		
		\begin{remark} \thlabel{remark:4areisom}
			Note that in the case when $N = 4$ the actions of $\Lambda_1$ and $\Lambda_3$ on $K[x_0, x_1, x_0', x_1']$ have the same invariants. We therefore treat both cases simultaneously.  
			
			There is also an explanation for this phenomenon in terms of the moduli interpretation of $Z(4, r)$. Every elliptic curve $E/K$ is $(4, 3)$-congruent with its quadratic twist by its discriminant (see \cite[Corollary~7.4]{BD_TAOG2CW44SJ} or \cite[Proposition~3.6]{F_COECAFNSMNGR} for a more general statement). In particular, if $E$ and $E'$ are $(4, 1)$-congruent, then the quadratic twist of $E$ by its discriminant $E^{\Delta(E)}$ is $(4, 3)$-congruent to $E'$. This induces an isomorphism $Z(4, 1) \cong Z(4, 3)$ which commutes with the $j$-invariant maps $Z(4, r) \to Z(1)$.
		\end{remark}

		Recall the forms $D = D(x_0, x_1)$, $c_4 = c_4(x_0, x_1)$, and $c_6 = c_6(x_0, x_1)$ defined in Section~\ref{sec:action}. We let $D' = D(x_0', x_1')$, $c_4' = c_4(x_0', x_1')$, and $c_6' = c_6(x_0', x_1')$. The forgetful map $X(N) \times X(N) \to Z(1)$ is given by $([x_0 : x_1], [x_0' : x_1']) \mapsto (j, j')$ where $j = c_4^3/D^N$ and $j' =c_4'^3/D'^N$. We write $J = j/1728$ (similarly $J' = j'/1728$).
	
		Additionally, we define the following $r$--bi-invariant forms:
		\begin{equation*}
			\begin{aligned}
				(N,r) = (2, 1) :& \qquad I_{1,1} = 192x_0 x_0' + x_1 x_1', \\
				(N,r) = (4, r) :& \qquad I_{2,2} = 16(x_0 x_1' - x_1 x_0')^2,\\
								& \qquad I_{4,4} = 256 x_0^4 x_0'^4 + 16 x_0^4 x_1'^4 + 192 x_0^2 x_1^2 x_0'^2 x_1'^2 + 16 x_1^4 x_0'^4 + x_1^4 x_1'^4,\\
				(N,r) = (3, 1) :& \qquad I_{2,2} = 3(54 x_0^2 x_0'^2 + x_0 x_1 x_1'^2 + x_1^2 x_0' x_1'),\\
							    & \qquad I_{6,6} = 3^6 (x_0x_1' - x_1x_0')^6 ,\\
				(N,r) = (3, 2) :& \qquad I_{2,2} = (18x_0 x_0' + x_1 x_1')^2.\\
			\end{aligned}
		\end{equation*}

		We then make the following choices of symmetric $r$--bi-invariants:

		\begin{table}[H]
		\centering 
		\begin{tabular}{c|c}
			$(N,r)$		& $r$--bi-invariants $(w_0, ..., w_k)$ \\
			\hline
			$(2,1)$		&	$\left( I_{1,1}^3 ,\enspace 2 I_{1,1} c_4 c_4' , \enspace 1728 D D' \right)$ \\
			$(4,r)$		&	$\left( I_{4,4}^3 ,\enspace 2 I_{4,4} c_4 c_4',\enspace 3 I_{2,2}^2 I_{4,4}^2,\enspace 72 I_{2,2} I_{4,4} D D' \right)$\\	
			$(3,1)$		&	$\left( 12I_{6,6}, \enspace 4 I_{2,2}^3, \enspace c_6 c_6',\enspace 144 I_{2,2} DD' \right) $ \\
			$(3,2)$		&	$\left(\frac{3}{2} (3 I_{2,2}^2 - c_4c_4' - 144 DD'), \enspace c_4 c_4', \enspace 144 D D'\right)$ \\
		\end{tabular}
    \caption{Choice of $r$--bi-invariants inducing a rational map $W(N, r) \dashrightarrow \bbP^k$.}
		\end{table}

		\subsubsection{Case \texorpdfstring{$(N, r) = (2, 1)$}{(N,r)=(2,1)}} \label{subsubsec:21} First note that $c_6 c_6' = 4w_0 - \frac{3}{2}w_1 - w_2$. In particular we have
		\begin{equation*}
			JJ' = \frac{w_1^3}{8 w_0 w_2^2},
		\end{equation*}
		\begin{equation*}
			(J - 1)(J' - 1) = \left( \frac{4w_0 - \frac{3}{2}w_1 - w_2}{w_2} \right)^2.
		\end{equation*}

		Let $S(2, 1) = \bbP^2$ with coordinates $w_0, w_1, w_2$. Then there is a commutative diagram 
		\[\begin{tikzcd}
			{W(2, 1)} & {S(2, 1)}\\
			{} & {W(1)}			
			\arrow[dashed, "\psi_{2,1}", from=1-1, to=1-2]
			\arrow[from=1-1, to=2-2]
			\arrow[dashed, from=1-2, to=2-2]
		\end{tikzcd}\] 
		where the map $S(2,1) \dashrightarrow W(1)$ is given by the rational functions $JJ'$ and $(J-1)(J'-1)$ above.

		\subsubsection{Case \texorpdfstring{$(N, r) = (4, r)$}{(N,r)=(4,r)}} \label{subsubsec:4r} We first note that the invariants $w_0, w_1, w_2, w_3$ satisfy the single relation
    \begin{equation*}
      192 w_0^2 - 96 w_0 w_1 + 128 w_0 w_3 + 48 w_0 w_2 - w_2^2 = 0.
    \end{equation*}
    Therefore the image of the map $\varphi_{4,r} \colon X(4) \times X(4) \dashrightarrow \bbP^3$ is the quadric surface $S(4,r) = \{ 192 w_0^2 - 96 w_0 w_1 + 128 w_0 w_3 + 48 w_0 w_2 - w_2^2 \} \subset \bbP^3$.
		
		Let $H = \ker \left( \SL_2(\bbZ/4\bbZ) \to \SL_2(\bbZ/2\bbZ) \right)$, so that the forgetful map $X(4) \to X(2)$ is given by taking the quotient by the action of $H/\{\pm 1\}$. Note that the $r$--bi-invariants $I_{4,4}$, $c_4 c_4'$, and $(DD')^2$ are invariant under the action of $H \times H$. In particular $w_0$, $w_1$ and $w_3^2/w_2$ are invariant under $H \times H$ and we have the relation $c_6 c_6' = 4w_0 - \frac{3}{2}w_1 - (w_3^2/w_2)$.

		Using the equations for the forgetful map $X(4) \to X(2)$ given in Section~\ref{sec:action} we see that the diagram
		\[\begin{tikzcd}
			{W(4, r)} & {S(4, r)}\\
			{W(2, 1)} & {S(2, 1)}			
			\arrow[dashed, "\psi_{4,r}", from=1-1, to=1-2]
			\arrow[from=1-1, to=2-1]
			\arrow[dashed, "\psi_{2,1}", from=2-1, to=2-2]
			\arrow[dashed, from=1-2, to=2-2]
		\end{tikzcd}\]
		commutes, where the map $S(4, r) \dashrightarrow S(2, 1)$ is given by $[w_0 : w_1 : w_2 : w_3] \mapsto [w_0 : w_1 : w_3^2/w_2]$.

		\subsubsection{Case \texorpdfstring{$(N, r) = (3, 1)$}{(N,r)=(3,1)}} \label{subsubsec:31} The image of the rational map $\varphi_{3,1} \colon X(3) \times X(3) \dashrightarrow \bbP^3$ is the cubic surface $S(3, 1) = \{8 w_0^2 w_1 - 60 w_0 w_1^2 + 12 w_0 w_1 w_2 + 36 w_0 w_1 w_3 - 9 w_1^3 + 27 w_1^2 w_3 - 27 w_1 w_3^2 + 9 w_3^3 = 0\}$. Additionally there are the relations:
		\begin{equation*}
			JJ' = \frac{w_1(4 w_0^2 - 192 w_0 w_1 + 12 w_0 w_2 + 72 w_0 w_3 + 603 w_1^2 - 144 w_1 w_2 - 918 w_1 w_3 + 9 w_2^2 + 108 w_2 w_3 + 351 w_3^2)}{36 w_3^3},
		\end{equation*}
		\begin{equation*}
			(J-1)(J'-1) = \frac{w_1 w_2^2}{4w_3^3}.			
		\end{equation*}
		We therefore have a commutative diagram
		\[\begin{tikzcd}
			{W(3, 1)} & {S(3,1)}\\
			{} & {W(1)}			
			\arrow[dashed, "\psi_{3,1}", from=1-1, to=1-2]
			\arrow[from=1-1, to=2-2]
			\arrow[dashed, from=1-2, to=2-2]
		\end{tikzcd}\] 
		where the morphism $S(3, 1) \dashrightarrow W(1)$ is given by the rational functions $JJ'$ and $(J-1)(J'-1)$ above.

		\subsubsection{Case \texorpdfstring{$(N, r) = (3, 2)$}{(N,r)=(3,2)}} \label{subsubsec:32} Let $S(3,2) = \bbP^2$. The bi-invariants $w_0$, $w_1$, and $w_2$ are linearly independent, so the map $\varphi_{3,2} \colon X(3) \times X(3) \dashrightarrow S(3,2)$ is surjective. We have the relations 
		\begin{equation*}
			JJ' = \left( \frac{w_1}{w_2} \right)^3 ,
		\end{equation*}
		\begin{equation*}
			(J-1)(J'-1) = \frac{(w_0^2 - 3 w_1^2 - 3 w_1 w_2 - 3 w_2^2)^2}{4 w_2^3 (2 w_0 + 3 w_1 + 3 w_2)}.
		\end{equation*}
		Therefore there is a commutative diagram
		\[\begin{tikzcd}
			{W(3, 2)} & S(3,2)\\
			{} & {W(1)}			
			\arrow[dashed, "\psi_{3,2}", from=1-1, to=1-2]
			\arrow[from=1-1, to=2-2]
			\arrow[dashed, from=1-2, to=2-2]
		\end{tikzcd}\] 
		where the map $S(3,2) \dashrightarrow W(1)$ is given by the rational functions $JJ'$ and $(J-1)(J'-1)$ above.

	\subsection{Symmetric conditions for \texorpdfstring{$2$, $3$, and $4$}{2,3, and 4}-congruences}

  We now prove \thref{thm:plys}. This provides a simple characterisation of $2$, $3$, and $4$-congruent elliptic curves in terms of their $j$-invariants. It is important for our computations of the surfaces $Z(12, r)$ that, unlike the equations of Fisher~\cite{F_THOAG1C} and of Rubin and Silverberg~\cite{RS_FOECWCMPR,RS_M2ROEC,S_EFOECWPMNR}, the conditions in \thref{thm:plys} are symmetrical in the pair of congruent elliptic curves, $E$ and $E'$. This will be useful since we will be able to work with the surfaces $W(12, r)$ which have simpler birational geometry. 
		
		The case when $N = 2$ appears in \cite[Lemma~3.5]{F_EMSFCOEC} with a different (computational) proof using the equations for $X_E^1(2)$. It is also possible to verify \thref{thm:plys} using the equations of Fisher and of Rubin and Silverberg for $X_E^r(N)$ (the necessary transformations are recorded in the electronic data~\cite{ME_ELECTRONIC_Z12}). We opt for a more geometric approach which explains how we found these conditions.
		
		\begin{theorem} \thlabel{thm:plys}
			Let $K$ be a field of characteristic not equal to $2$ or $3$. Let $E/K$ and $E'/K$ be elliptic curves with $j$-invariants $j$ and $j'$ respectively. Let $c_4, c_6, \Delta$ and $c_4', c_6', \Delta'$ be the $c$-invariants and discriminant of (a Weierstrass model of) $E$ and $E'$ respectively. Suppose that $j, j' \not\in \{0,1728\}$. Let $J = j/1728$ and $J' = j'/1728$. Then we have the following:
			\begin{enumerate}[label=\eniii]
				\item \label{thm:plys_2-1} $E$ and $E'$ are $2$-congruent if and only if there exist $\alpha, \beta \in K$ such that 
				\begin{equation*}	\label{eqn:2cong_1}
					{\alpha}^2 - (J -1)(J'-1) = 0,
				\end{equation*}
				\begin{equation*}	\label{eqn:2cong_2}
					{\beta}^3 - 3 J J' {\beta} - 2 J J' ({\alpha}+1) =0.
				\end{equation*}
				
				\item \label{thm:plys_4-r} If $J \neq J'$ then $E$ and $E'$ are $(4, r)$-congruent if and only if there exist $\alpha, \beta, \gamma, \tau \in K$ such that $\alpha$ and $\beta$ satisfy the conditions of \ref{thm:plys_2-1} and 
				\begin{equation*}	\label{eqn:4cong_1}
					\gamma^4 - 6 J J' \beta \gamma^2 - 16 (J J')^2 \gamma + 3 (J J')^2 (4 J J' -  \beta^2) = 0,
				\end{equation*}
				\begin{equation*}	\label{eqn:4cong_2}
					\tau^2 = \begin{cases}
						3 c_6 c_6' \delta \alpha \qquad &\textnormal{if } r \equiv 1 \pmod{4}, \\
						3 c_6 c_6' \alpha \qquad &\textnormal{if } r \equiv 3 \pmod{4}
					\end{cases}
				\end{equation*}
				where $\delta = 3(JJ' - (\alpha + 1)^2)$.
				
				\item \label{thm:plys_3-1} If $J \neq J'$ and $(\sqrt[3]{J} + 1)(\sqrt[3]{J'} + 1) \neq 1$ (for all choices of cube roots) then $E$ and $E'$ are $(3, 1)$-congruent if and only if there exist $\alpha, \beta, \tau \in K$ such that 
				\begin{equation*}	
					\alpha^3 - 3  J J' \alpha -  J J' ( J + J')=0,
				\end{equation*}
				\begin{equation*}	
					\beta^4 - 6 (J - 1)(J' - 1) \beta^2 - 8 ( J - 1)^2 (J' - 1)^2 \beta - 3 (4\alpha + 1)( J - 1)^2 (J' - 1)^2=0,
				\end{equation*}
				\begin{equation*}	
					\tau^2 = 3 c_6 c_6' \delta,
				\end{equation*}
				where 
				\begin{equation*}
					\delta = -6\left(\frac{2\beta^3 - (5\alpha + 2)\beta^2 - 10(J-1)(J'-1)\beta + 3(J-1)(J'-1)(13\alpha - 2 + 6(J + J'))}{\beta^3 - 3(J-1)(J'-1)\beta - 2(J-1)^2(J'-1)^2}\right).
				\end{equation*}			
				
				\item \label{thm:plys_3-2}If $J \neq J'$ then $E$ and $E'$ are $(3, 2)$-congruent if and only if there exists $\alpha, \beta, \tau \in K$ such that 
				\begin{equation*}
					\alpha^3 -  J J'=0, 
				\end{equation*}
				\begin{equation*}	
					\beta^4 -6 (\alpha + 1) ( J - 1) (J' - 1) \beta^2 -8 ( J - 1)^2 (J' - 1)^2  \beta - 3 (\alpha - 1)^2 ( J - 1)^2 (J' - 1)^2=0 , 
				\end{equation*}
				\begin{equation*}	
					\tau^2 = \begin{cases}
						3 c_6 c_6' \beta 	\qquad &\textnormal{if } \beta \neq 0, \\
						-2 c_6 c_6' 		\qquad &\textnormal{if } \beta = 0.
					\end{cases}
				\end{equation*}
			\end{enumerate}
			
		\end{theorem}
		
		\begin{remark} Note that in \thref{thm:plys}:
			
			\begin{enumerate}[label=\eniii]
				\item The conditions on $\alpha$ are convenient ways of writing the conditions that $\Delta \Delta'$ is a square when $N = 2, 4$ and that $\Delta\Delta'$ is a cube when $(N, r) = (3, 2)$. In the case when $(N,r) = (3,1)$ if we define $\alpha' = \frac{1}{J - J'} \left( \alpha^2/J' - \alpha - 2J\right)$ then $(\alpha')^3 = J/J'$. In particular when $(N,r) = (3,1)$ the condition on $\alpha$ is equivalent to $\Delta/\Delta'$ being a cube.
				
				When $N = 4$ the term $\delta$ is equal to $\Delta$ (and $\Delta'$) up to a square factor (modulo squares we have $\delta/\Delta \equiv 3\delta(J - 1) \equiv (\alpha - (J - 1))^2$). The conditions on $\tau$ are therefore equivalent to the fact that $\Delta/\Delta'$ is a fourth power when $(N, r) = (4, 1)$ and  $\Delta\Delta'$ is a fourth power when $(N, r) = (4, 3)$.
				
				\item When $K$ contains an $N^\text{th}$ root of unity, it is clear that the polynomials in $\alpha$ in \thref{thm:plys} may have more than one root in $K$. Indeed, the theorem would be false if we required the existence of a $\beta$ (and $\gamma, \tau$) \emph{for a fixed choice of $\alpha$ (and $\beta$)}. For example, the curves with LMFDB labels \LMFDBLabel{196.b1} and \LMFDBLabel{196.b2} are $3$-isogenous and therefore $(4,3)$-congruent. In this example the polynomial in $\beta$ has no root when the negative root $\alpha$ is chosen. When $\alpha$ is chosen to be positive the polynomial in $\beta$ splits, but the polynomial in $\gamma$ has a root for only one of the possible choices of $\beta$.
				
				When $(N, r) = (3, 1)$ a similar example is given over $\bbQ(\zeta_3)$ by the pair \LMFDBLabel{15.a1} and \LMFDBLabel{15.a5} and when $(N, r) = (3, 2)$ by the pair \LMFDBLabel{11.a1} and \LMFDBLabel{11.a2}.				
				
				\item When $(N, r) = (3, 1)$ the condition that $(\sqrt[3]{J} + 1)(\sqrt[3]{J'} + 1) \neq 1$ is necessary. The elliptic curves with LMFDB labels \LMFDBLabel{245.a1} and \LMFDBLabel{1323.s1} are not (quadratic twists of) $(3, 1)$-congruent elliptic curves, however there do exist $\alpha, \beta \in \bbQ$ satisfying \threfpart{thm:plys}{3-1}. Moreover when $N = 3, 4$ the theorem is false when $J = J'$. Counterexamples are given by the pairs of elliptic curves \LMFDBLabel{14.a1} and \LMFDBLabel{126.b1} when $(N,r) = (3,1)$, \LMFDBLabel{245.a1} and \LMFDBLabel{245.b1} when $(N,r) = (3,2)$, \LMFDBLabel{14.a6} and \LMFDBLabel{112.c6} when $(N, r) = (4, 1)$, and \LMFDBLabel{14.a6} and \LMFDBLabel{784.b6} when $(N, r) = (4, 3)$.
			\end{enumerate}
		\end{remark}
		
		\begin{proof}[Proof of \thref{thm:plys}]
			If $j = j'$ then $E$ and $E'$ are $2$-congruent (since an elliptic curve is $2$-congruent to any quadratic twist and $j, j' \not \in \{0, 1728\}$ by assumption). Condition (i) is always satisfied in this case, so we are free to assume that $j \neq j'$ when $(N, r) = (2, 1)$.
			
			For each $(N, r)$ let $P_{\alpha} \in K[J, J', \alpha]$, $P_{\beta} \in K[J, J', \alpha, \beta]$ (and $P_{\gamma} \in K[J, J', \alpha, \beta, \gamma]$ when $N = 4$) be the polynomials from the statement of the theorem. 
			
			When $N = 2, 3$ let $\mathcal{Z}(N, r) = \Spec K[J, J', \alpha, \beta]/(P_{\alpha}, P_{\beta})$ and $\mathcal{Z}(4, r) = \Spec K[J, J', \alpha, \beta, \gamma]/(P_{\alpha}, P_{\beta}, P_{\gamma})$. Similarly when $N = 2, 3$ let $\mathcal{W}(N, r) = \Spec K[JJ', (J-1)(J' -1), \alpha, \beta]/(P_{\alpha}, P_{\beta})$ and let $\mathcal{W}(4, r) = \Spec K[JJ', (J-1)(J' -1), \alpha, \beta, \gamma]/(P_{\alpha}, P_{\beta}, P_{\gamma})$.
			
			Write, in the notation of Section~\ref{sec:biinvs}:
			\begin{align*}
				(N,r) &= (2, 1) :&	\alpha &= \frac{4w_0 - \frac{3}{2}w_1 - w_2}{w_2} & 	\beta &= w_1/w_2 \\
				(N,r) &= (4, r) :&	\alpha &= \frac{4w_0 - \frac{3}{2}w_1 - (w_3^2/w_2)}{(w_3^2/w_2)} & 	\beta &= \frac{w_1}{w_3^2/w_2}	& \gamma = \frac{JJ'w_2}{\beta w_3}  \\
				(N,r) &= (3, 1) :&	\alpha &=  \frac{\eta}{12 w_0 w_3^2} & 	\beta &= \frac{(\frac{3}{2}(3w_1 - w_3)/w_3)^2 - 3(\alpha + 1)}{2} \\					
				(N,r) &= (3, 2) :&	\alpha &=  w_1/w_2 & 	\beta &= \frac{ \left(w_0/w_2\right)^2 - 3( \alpha^2 + \alpha + 1)}{2} \\
			\end{align*}
			where $\eta = (w_0 + 3 w_1 - 3 w_3)(8 w_0 w_1 - 3 w_1^2 + 6 w_1 w_3 - 3 w_3^2)$. Then these choices induce maps $S(N, r) \dashrightarrow \mathcal{W}(N, r)$ (where $S(N,r)$ are the surfaces defined in Sections~\ref{subsubsec:21}--\ref{subsubsec:32}) and therefore maps $W(N,r) \to \mathcal{W}(N,r)$. 
			
			In particular we have a commutative diagram  
			\[\begin{tikzcd}
				Z(N,r) && \mathcal{Z}(N,r) \\
				& Z(1)
				\arrow["\mathcal{J}_{Z}"', from=1-1, to=2-2]
				\arrow["\mathcal{J}_{\mathcal{Z}}", from=1-3, to=2-2]
				\arrow["\pi", dashed, from=1-1, to=1-3]
			\end{tikzcd}\]
			where $\mathcal{J}_{Z}$ and $\mathcal{J}_{\mathcal{Z}}$ are morphisms. The degree of the morphism $\mathcal{J}_Z$ is equal to $\lvert \SL_2(\bbZ/N\bbZ) \rvert / 2$ and the degree of $\mathcal{J}_{\mathcal{Z}}$ may be read off from the degrees of the polynomials $P_\alpha$, $P_\beta$, and $P_\gamma$. Since these degrees are equal for each $(N, r)$ the map $\pi$ is birational. 
		
			Write $R = K[J, J']\left[\frac{1}{JJ'(J-1)(J'-1)(J-J') \xi }\right]$ where 
			\begin{equation*}
				\xi = 	\begin{cases}
									1 & \text{if $N \in \{2, 4\}$, or $(N, r) = (3, 2)$} \\
									(J J')^3 + 3(J J')^2 (J + J') - 27 (J J')^2 + 3 (J J') (J + J')^2 + (J + J')^3 & \text{if $(N, r) = (3, 1)$.} \\
								\end{cases}
			\end{equation*}
			Letting $Z'(1) = \Spec R$ we have a natural open immersion $Z'(1) \subset Z(1)$. Let $Z'(N, r)$ and $\mathcal{Z}'(N, r)$ be the preimages of $Z'(1)$ under the maps $\mathcal{J}_{Z}$ and $\mathcal{J}_{\mathcal{Z}}$ respectively. 
			
			Explicitly $Z'(N,r)$ and $\mathcal{Z}'(N, r)$ are the surfaces given by deleting the points on $Z(N,r)$ and $\mathcal{Z}(N,r)$ above the loci:
			\begin{align*}
				(N,r) &= (2, 1) :&	&\text{$J = J'$ and $J, J' \in \{0, 1\}$,} \\
				(N,r) &= (4, r) :&	&\text{$J = J'$ and $J, J' \in \{0, 1\}$,}  \\
				(N,r) &= (3, 1) :&	& \text{$J = J'$, $J, J' \in \{0, 1\}$, and $(\sqrt[3]{J} + 1)(\sqrt[3]{J'} + 1) = 1$,} \\
				(N,r) &= (3, 2) :&	&\text{$J = J'$ and $J, J' \in \{0, 1\}$.} \\
			\end{align*}
			
			We next show that $\pi$ extends to an isomorphism $\widetilde{\pi} \colon Z'(N,r) \to \mathcal{Z}'(N, r)$.  Note that the morphism $Z'(N, r) \to Z'(1)$ is finite \'{e}tale since $X(N) \to X(1)$ is finite \'{e}tale away from $j = 0 , 1728, \infty$. 
			
			When $N = 2, 3$ we see that the finite morphism $\mathcal{Z}'(N, r) \to Z'(1)$ factors via 
			\[\begin{tikzcd}
				\mathcal{Z}'(N, r) = \Spec {R[\alpha, \beta]/(P_\alpha, P_\beta)} & \Spec {R[\alpha]/(P_\alpha)} & Z'(1).
				\arrow[from=1-1, to=1-2]
				\arrow[from=1-2, to=1-3]
			\end{tikzcd}\]
			and similarly when $N = 4$ with an extra factor accounting for the polynomial in $\gamma$. In particular when $N = 2, 3$, to show that $\mathcal{Z}'(N, r) \to Z'(1)$ is \'{e}tale it suffices to show that $\partial P_\alpha / \partial \alpha$ and $\partial P_\beta / \partial \beta$ are units in $R[\alpha]/(P_\alpha)$ and $R[\alpha, \beta]/(P_\alpha, P_{\beta})$ respectively. When $N = 4$ we must also check that $\partial P_\gamma / \partial \gamma$ is a unit in $R[\alpha, \beta, \gamma]/(P_\alpha, P_{\beta}, P_{\gamma})$.

			Computing successive resultants in \texttt{Magma} and using the fact that $2$, $3$, $J$, $J'$, $J-1$, $J'-1$, $J - J'$, and $\xi$ are units in $R$ shows that the ideals $\left( P_{\alpha}, \frac{\partial P_\alpha}{\partial \alpha} \right)$, $\left( P_{\alpha}, P_\beta, \frac{\partial P_\beta}{\partial \beta} \right)$, and $\left( P_{\alpha}, P_\beta, P_\gamma, \frac{\partial P_\gamma}{\partial \gamma} \right)$ contain $1$. In particular, the morphisms $\mathcal{Z}'(N, r) \to Z'(1)$ are finite \'{e}tale. 
			
			Applying the following lemma with $X = \mathcal{Z}'(N, r)$, $Y = Z'(N, r)$, and $S = Z'(1)$ shows that $\pi$ extends to an isomorphism $\widetilde{\pi} \colon \mathcal{Z}'(N, r) \to Z'(N, r)$.

			\begin{lemma} \thlabel{lemma:birtoisom}
				Let $X$, $Y$, and $S$ be integral $K$-varieties. Suppose that $S$ is smooth and that there exists a commutative diagram 
				\[\begin{tikzcd}
					X && Y \\
					& S
					\arrow["p_X"', from=1-1, to=2-2]
					\arrow["p_Y", from=1-3, to=2-2]
					\arrow["\pi", dashed, from=1-1, to=1-3]
				\end{tikzcd}\]
				where $p_X$ and $p_Y$ are finite \'{e}tale morphisms and $\pi$ is birational. Then $\pi$ extends uniquely to an isomorphism $\widetilde{\pi} \colon X \to Y$.

				\begin{proof}
					Because $\pi$ is birational there exist dense open subschemes $U \subset X$ and $V \subset Y$ such that $\pi \vert_U \colon U \to V$ is an isomorphism. Since $p_X$ and $p_Y$ are finite \'{e}tale we may choose $U$ and $V$ such that $p_X \vert_U \colon U \to S$ and $p_X \vert_V \colon V \to S$ are finite \'{e}tale.

					By assumption $S$ is integral and nonsingular, so the natural functor from the category of finite \'{e}tale covers of $S$ to the category of finite \'{e}tale covers of $p_X(U)$ is fully faithful (see \cite[\href{https://stacks.math.columbia.edu/tag/0BQG}{Lemma 0BQG}]{stacks-project} noting the proof of \cite[\href{https://stacks.math.columbia.edu/tag/0BQI}{Lemma 0BQI}]{stacks-project}). In particular the isomorphism $\pi\vert_U \colon U \to V$ extends uniquely to a morphism $\widetilde{\pi} \colon X \to Y$. A symmetrical argument shows that the inverse of $\pi \vert_U$ extends uniquely to a morphism $Y \to X$. Therefore by uniqueness $\widetilde{\pi}$ is an isomorphism.
				\end{proof}
			\end{lemma}

			The surface $Z(N, r)/K$ acts as a fine moduli space for pairs of elliptic curves $E/K$ and $E'/K$ with $j(E), j(E') \not\in \{0, 1728\}$ which are $(N, r)$-congruent up to a simultaneous quadratic twist (or simply up to quadratic twist if $N = 2$). Therefore we have only shown that elliptic curves $E/K$ and $E'/K$ satisfying the conditions on $\alpha$, $\beta$ (and $\gamma$ is $N = 4$) are $(N, r)$-congruent \emph{up to quadratic twist}. When $N = 3, 4$ it remains to show that the condition on $\tau$ is satisfied if and only if $E$ and $E'$ are $(N,r)$-congruent.
			
			The following lemma will allow us to identify the correct quadratic twist.
			
			\begin{lemma}[cf., {\cite[Proposition~13(A)]{FK_OTSTOIOTPTOEC}}] \thlabel{lemma:badfactors}
				Let $N \geq 3 $ be an integer and let $K$ be the field of fractions of a complete discrete valuation ring with residue field of characteristic coprime to $N$. Suppose that $E/K$ has good reduction and $E'/K$ has potential good reduction. If $E$ and $E'$ are $N$-congruent over $K$ then $E'$ has good reduction over $K$.
				
				\begin{proof}
					Let $K^{\text{ur}}$ be the maximal unramified extension of $K$ and let $L = K^{\text{ur}}(E[N])$ and $L' = K^{\text{ur}}(E'[N])$. Since $E$ and $E'$ are $N$-congruent $L = L'$.

					But $L$ (resp. $L'$) is the smallest extension of $K^{\text{ur}}$ over which $E$ (resp. $E'$) obtains good reduction (see~\cite[\S2 Corollary~3]{ST_GROAV}). Since $E$ has good reduction over $K$ we have $L' = L = K^{\text{ur}}$. Therefore $E'$ has good reduction over $K$. 
				\end{proof}
			\end{lemma}

			The fibres of the maps $\mathcal{Z}'(N, r) \xrightarrow{j} X(1)$ have genus $0$ and are birational to the curves $X_E^r(N)$ where $E$ is an elliptic curve with $j$-invariant $j$. By parametrising the generic fibre we obtain a birational map $\lambda \colon \mathcal{Z}'(N, r) \to X(1) \times \bbP^1$ which is an isomorphism onto its image since the fibres of $\mathcal{Z}'(N, r) \xrightarrow{j} X(1)$ are nonsingular.

			Let $\mathscr{E} \to \mathcal{Z}'(N, r)$ be given by the Weierstrass equation
			\begin{equation} \label{eqn:scrE}
				y^2 = x^3 - 27 \frac{j}{j-1728}x - 54 \frac{j}{j-1728} 
			\end{equation}
			and let $\mathscr{E}' \to \mathcal{Z}'(N, r)$ be given by the Weierstrass equation
			\begin{equation} \label{eqn:scrE'}
				y^2 = x^3 - 27d^2 \frac{j'}{j'-1728} x - 54d^3 \frac{j'}{j'-1728}.
			\end{equation}
			where 
			\begin{align*}
				(N,r) &= (4, 1) :&	d &= \frac{3 \delta \alpha j j'}{(j - 1728)(j' - 1728)}, \\
				(N,r) &= (4, 3) :&	d &= \frac{3 \alpha j j'}{(j - 1728)(j' - 1728)},  \\
				(N,r) &= (3, 1) :&	d &= \frac{3 \delta j j'}{(j - 1728)(j' - 1728)}, \\					
				(N,r) &= (3, 2) :&	d &= \frac{3 \beta j j'}{(j - 1728)(j' - 1728)}. \\
			\end{align*}
			
			By construction the rational functions $j$, $j'$, $(j-1728)$, and $(j'-1728)$ are units in $R$. We check in \texttt{Magma} that when $N = 4$ the functions $\delta$ and $\alpha$ are units in $R[\alpha, \beta, \gamma]/(P_\alpha, P_{\beta}, P_{\gamma})$. Similarly when $(N, r) = (3, 1)$ the function $\delta$ is a unit in $R[\alpha, \beta]/(P_\alpha, P_{\beta})$. Computing the discriminants of the Weierstrass equations in (\ref{eqn:scrE}) and (\ref{eqn:scrE'}) we see that the fibres of $\mathscr{E}, \mathscr{E}' \to \mathcal{Z}'(N, r)$ are nonsingular except when $(N, r) = (3, 2)$ and $\beta = 0$.

			The generic fibres of $\mathscr{E}, \mathscr{E}' \to \mathcal{Z}'(N, r)$ are elliptic curves $\mathcal{E}/\bbQ(j,t)$ and $\mathcal{E}'/\bbQ(j,t)$ and we may choose models such that $\Delta(\mathcal{E}), \Delta(\mathcal{E}') \in \bbZ[j,t]$, where $j$ and $t$ are coordinates on $X(1) \times \bbP^1$. 
			
			Since $\mathcal{E}$ and $\mathcal{E}'$ are $(N, r)$-congruent up to quadratic twist there exists some squarefree $D \in \bbZ[j,t]$ such $\mathcal{E}$ is $(N, r)$-congruent to $(\mathcal{E}')^D$. If $D$ has an irreducible factor $q$ not dividing $\Delta(\mathcal{E})\Delta(\mathcal{E}')$ then applying \thref{lemma:badfactors} over the completion of $\bbZ[j,t]$ at $q$ we see that $\mathcal{E}$ and $(\mathcal{E}')^D$ are not $N$-congruent, hence $D$ divides $\Delta(\mathcal{E})\Delta(\mathcal{E}')$. In particular, there are finitely many possibilities for $D \in \bbZ[j,t]$. 
			
			Explicitly computing the parametrisation $\lambda$ in \texttt{Magma}, specialising at a small number of $j,t \in \bbQ$, and comparing traces of Frobenius is enough to prove that $D = 1$ when $N = 3$, and $D = 1$ or $\delta$ when $N = 4$. 
			
			The condition on $\tau$ follows immediately when $(N, r) = (3, 1)$ and when $(N, r) = (3, 2)$ and $\beta \neq 0$. Moreover we see that elliptic curves are $4$-congruent if and only if they satisfy one of the conditions on $\tau$ in \ref{thm:plys_4-r}. 
			
			Note that $3 c_6 c_6' \alpha$ is not a square for the pair of $5$-isogenous  (hence $(4, 1)$-congruent) elliptic curves with LMFDB labels \LMFDBLabel{38.b1} and \LMFDBLabel{38.b2}. Therefore $D = 1$ when $(N, r) = (4, 1)$. Similarly $3 c_6 c_6' \delta \alpha$ is not a square for the pair of $3$-isogenous elliptic curves \LMFDBLabel{44.a1} and \LMFDBLabel{44.a2}, so $D = 1$ when $(N, r) = (4, 3)$.

			The curve on $\mathcal{Z}'(3, 2)$ given by the vanishing of $\beta$ admits an obvious parametrisation by setting $J = t$, $J' = 1/t$, $\alpha = 1$, and $\beta = 0$. Repeating the above calculation with $d = \frac{-2jj'}{(j-1728)(j'-1728)}$ completes the proof.
		\end{proof}

		\begin{remark}
			For $N \leq 4$, rather than recovering $X_E^r(N)$ it is also possible to give a parametrisation $\bbA^2 \dashrightarrow Z(N, r)$ so that the involution swapping the roles of $E$ and $E'$ is given by swapping the coordinates on $\bbA^2$. In the electronic data~\cite{ME_ELECTRONIC_Z12} we record these symmetric families of $(N, r)$-congruent elliptic curves over $\bbQ(a,b)$. These parametrisations were found by first parametrising the surfaces $S(N, r)$ which are birational to $W(N,r)$.
		\end{remark}

\section{Constructing \texorpdfstring{$Z(N,r)$}{Z(N,r)} via fibre products} \label{sec:fibprod}
	
	Let $N_1, N_2 \geq 2$ be coprime integers and let $N = N_1 N_2$. We define $Z^+(N, r)$ to be the fibre product  
	\[\begin{tikzcd}
		& {Z^+(N, r)} \\
		{Z(N_1, r)} && {Z(N_2, r)} \\
		& {Z(1)}
		\arrow[from=1-2, to=2-1]
		\arrow[from=1-2, to=2-3]
		\arrow[from=2-3, to=3-2]
		\arrow[from=2-1, to=3-2]
	\end{tikzcd}\]
	and similarly for $W^+(N, r)$.
	
	\begin{remark}
		Note that the surfaces $Z^+(N, r)$ and $W^+(N, r)$ may depend implicitly on the factorisation $N = N_1 N_2$.
	\end{remark}
	
	The surface $Z^+(N, r)$ has the following moduli interpretation: for each point $P \in Z^+(N, r)(K)$ with $j(P), j'(P) \not\in \{0, 1728\}$ there exist pairs of elliptic curves $(E_1, E_1')$ and $(E_2, E_2')$ defined over $K$ such that for $i= 1,2$ we have 
	\begin{itemize}
		\item $j(P) = j(E_i)$ and $j'(P) = j(E_i')$, and
		\item $E_i$ and $E_i'$ are $(N_i, r)$-congruent over $K$.
	\end{itemize}
	
	The surface $W^+(N, r)$ is then obtained from $Z^+(N, r)$ by taking the quotient by the involution swapping the roles of $E_i$ and $E_i'$. 
	
	\begin{remark} \thlabel{rmk:4isom}
		Suppose that $N_1 = 4$ and that $r_1, r_2 \in (\bbZ/N\bbZ)^\times$ are equal (up to a square) modulo $N_2$. Then the surfaces $Z^+(N, r_1)$ and $Z^+(N, r_2)$ are isomorphic since $Z(4, 1)$ and $Z(4, 3)$ are isomorphic (see \thref{remark:4areisom}) and the isomorphism respects the maps $Z(4, r) \to Z(1)$.
	\end{remark}
	
	By taking simultaneous quadratic twists of $E_2$ and $E_2'$ we may assume that $E_1 = E_2$. If $N_1 = 2$ then we immediately see that $Z^+(N, r)$ is birational to $Z(N, r)$ -- every elliptic curve is $2$-congruent to any quadratic twist, so we may replace $E_1'$ by its quadratic twist $E_2'$ so that $E_1$ is $(N, r)$-congruent to $E_1'$.
	
	Otherwise, assume that $N_i \neq 2$ for $i = 1,2$. By the above discussion if $\mathcal{E}/\bbQ(Z^+(N, r))$ is an elliptic curve with $j$-invariant $j$ then there exists a pair of elliptic curves $\mathcal{E}_1'/\bbQ(Z^+(N, r))$ and $\mathcal{E}_2'/\bbQ(Z^+(N, r))$ with $j$-invariant $j'$ such that $\mathcal{E}$ and $\mathcal{E}_i'$ are $(N_i, r)$-congruent. Let $d \in \bbQ(Z(N, r))$ be such that $\mathcal{E}_1'$ is isomorphic to the quadratic twist of $\mathcal{E}_2'$ by $d$. Then the modular diagonal quotient surface $Z(N, r)$ is birational to the double cover of $Z^+(N, r)$ whose function field is given by adjoining the square root of $d$ to $\bbQ(Z(N, r))$.
	
\section{Computing \texorpdfstring{$Z(12,r)$}{Z(12,r)}} \label{sec:compute}
	We now prove \thref{thm:Zequations} by applying \thref{thm:plys} to compute the fibre products described in Section~\ref{sec:fibprod}. To compute the surfaces $Z(12, r)$ for $r = 1, 7$ we use \threfpart{thm:plys}{2-1}--\ref{thm:plys_4-r} to construct the covers $Z(12, r) \to Z(3, 1)$. In contrast, to compute the surfaces $Z(12, r)$ for $r = 5, 11$ we use \threfpart{thm:plys}{3-2} to construct the covers $Z(12, r) \to Z(4, r)$. 
	
	\begin{remark}
		It is also possible to compute $Z(12, r)$ for $r = 1, 7$ by using \threfpart{thm:plys}{3-2} to construct the covers $Z(12, r) \to Z(4, r)$, and similarly when $r = 5, 11$ it is possible to compute $Z(12, r)$ by using \threfpart{thm:plys}{2-1}--\ref{thm:plys_4-r} to construct the covers $Z(12, r) \to Z(3, 2)$.
	\end{remark}

	\subsection{Choosing parametrisations via Cremona transformations}
		Since a model for $W(N,r)$ consists of both a birational model for the surface \emph{and} equations for the forgetful maps we naturally want to simplify both the equations for the surface and the maps. When the surface in question is rational we exploit the birational automorphisms of $\bbP^2$ to simplify the forgetful maps.
		
		It is well known that $\Bir_\bbQ(\bbP^2)$, the birational automorphism group of $\bbP^2$ over $\bbQ$, contains $\PGL_3(\bbQ)$ and the \emph{Cremona transformations} -- that is, birational automorphisms of the form
		\begin{equation*}
			[x_0:x_1:x_2] \mapsto [Q_0(x_0, x_1, x_2) : Q_1(x_0, x_1, x_2) : Q_2(x_0, x_1, x_2) ]
		\end{equation*}
		where $Q_0, Q_1, Q_2 \in \bbQ[x_0, x_1, x_2]$ are homogeneous forms of degree $2$ which vanish simultaneously at three $G_{\bbQ}$-conjugate points $q_1, q_2, q_3 \in \bbP^2(\overbar{\bbQ})$ in general position.
		
		Geometrically this operation blows up $q_1, q_2, q_3$ and blows down the three lines which contain exactly two of $q_1, q_2, q_3$.
		
		This suggests the following approach for simplifying a model for $W(N, r)$. Let $C/\bbQ$ be a curve in $\bbP^2$ given by the vanishing of a factor of the numerator of $JJ'$, $(J - 1)(J' - 1)$, or $(J - J')^2$. Applying a generic Cremona transformation to $\bbP^2$ will double the degree of $JJ'$, $(J - J')^2$, and $(J - 1)(J' - 1)$. If instead we apply a Cremona transformation centred at a $G_{\bbQ}$-stable triple of singular points $q_1, q_2, q_3 \in C(\Qbar)$ we may hope to decrease their degrees. 
		
		In the discussion that follows we often make choices of parametrisation which give rise to models for $W(N, r)$ with $j$-maps of small degree. In general, the parametrisations we give were not our first choice. Once we have a birational map $\bbA^2 \dashrightarrow W(N,r)$ (e.g., by spotting a genus $0$ fibration over $\bbP^1$), we apply iterative Cremona transformations as described above to simplify the $j$-maps.
	
	\subsection{Models for \texorpdfstring{$Z(12, 1)$}{Z(12,1)} and \texorpdfstring{$Z(12, 7)$}{Z(12,7)}} \label{sec:Z12_1_7}
		We will use \thref{thm:plys} to compute a model for the fibre product $W^+(12, r)$. 
		
		As in the proof of \thref{thm:plys} let $\mathcal{W}(3, 1)$ be the surface in $\bbA^4$ defined by the polynomials in \threfpart{thm:plys}{3-1}. We check in \texttt{Magma} that the rational map $\bbA^2 \dashrightarrow \mathcal{W}(3, 1)$ given by taking
		\begin{gather*}
			JJ' = \frac{-(-2 a b + b^{2} - 4 a - 4 b)^{3}}{(2 b + 1)^{3} (-a b + b^{2} - 2 a - 2 b)^{2}}, \\
			(J - 1)(J' - 1) = \frac{(b + 2)^{2} (a^{2} b - 2 a b^{2} + b^{3} + 2 a^{2} + a b - 2 b^{2} + a + b)^{2}}{(2 b + 1)^{3} (-a b + b^{2} - 2 a - 2 b)^{2}},
		\end{gather*}
		(and solving for $\alpha$ and $\beta$) is birational. In particular $W(3, 1)$ is birational to $\bbA^2$ with the forgetful map $W(3,1) \to W(1)$ given by the formulae for $JJ'$ and $(J - 1)(J' - 1)$ above.
		
		Let $\Delta = \Delta(\mathcal{E})$ and $\Delta' = \Delta(\mathcal{E}')$ be the discriminants of the elliptic curves $\mathcal{E}/\bbQ(Z(3,1))$ and $\mathcal{E}'/\bbQ(Z(3,1))$ with $j$-invariants $1728J$ and $1728J'$ respectively. We write $W(3, 1)^{\sqrt{\Delta\Delta'}}$ for the double cover of $W(3, 1)$ whose function field is given by adjoining $\sqrt{\Delta \Delta'}$ to $\bbQ(W(3,1))$.
		
		Then $W(3, 1)^{\sqrt{\Delta\Delta'}}$ is birational to the surface given by
		\begin{equation*}
			\alpha_4^2 - (J - 1)(J' - 1) = 0
		\end{equation*}
		in $\bbA^3$. By setting
		\begin{equation*}
			\alpha_4 = \frac{-(b + 2) (a^{2} b - 2 a b^{2} + b^{3} + 2 a^{2} + a b - 2 b^{2} + a + b)}{(2 b + 1)^{2} (-a b + b^{2} - 2 a - 2 b)} \alpha_4'
		\end{equation*}
		we see that $W(3, 1)^{\sqrt{\Delta\Delta'}}$ is birational to the surface given by 
		\begin{equation*}
			(\alpha_4')^2 - (2b + 1) = 0
		\end{equation*}
		in $\bbA^3$. We parametrise this surface by setting $(a,b,\alpha_4') = \left( c, \frac{d^2 - 1}{2}, d \right)$. Therefore $W(3, 1)^{\sqrt{\Delta\Delta'}}$ is birational to $\bbA^2$ and the forgetful map $W(3, 1)^{\sqrt{\Delta\Delta'}} \to W(3,1)$ is given by $(c, d) \mapsto \left(c, \frac{d^2 - 1}{2} \right)$.
		
		By \threfpart{thm:plys}{2-1} we may write $W(6, 1)$ as a degree $3$ cover of $\bbA^2$ given by the vanishing of 
		\begin{equation*}
			\beta_4^3 - 3 J J' \beta_4 - 2 J J' (\alpha_4 + 1) = 0
		\end{equation*}
		in $\bbA^3$. After making the change of coordinates 
		\begin{equation*}
			\beta_4 = \frac{ (d^{4} - 4 c d^{2} - 10 d^{2} - 12 c + 9)\left( ( d^{2} - 2 d - 3) \beta_4' - d^{4} + 2 c d^{2} + 2 d^{3} + 6 d^{2} + 6 c - 6 d - 9 \right)}{2 d^3 (d^{4} - 2 c d^{2} - 6 d^{2} - 6 c + 5) } \cdot \frac{1}{\beta_4'}
		\end{equation*}
		we find that $W(6, 1)$ is birational to the surface given by 
		\begin{equation*}
			(\beta_4')^{3} + 6 (\beta_4')^{2} - 3 (d + 1)(d - 3) \beta_4' - (2 c d^2 + 6 c - d^4 + 2 d^3 + 6 d^2 - 6 d - 9) = 0
		\end{equation*}
		in $\bbA^3$. We parametrise the above surface by setting
		\begin{equation*}
			(c, d, \beta_4') = \left(\frac{p^3 q + 6 p^2 q^2 + 9 p q^3 + 12 p q^2 - 12 p q + 9 q^4 + 12 q^3 - 24 q^2 - 16 q + 16}{2 q^{2} (3 q^2 + 4)}, \frac{2}{q}, \frac{p}{q}	\right). 
		\end{equation*}
		Hence $W(6, 1)$ is birational to $\bbA^2$ and the forgetful map $W(6, 1) \to W(3, 1)^{\sqrt{\Delta\Delta'}}$ is given by $(p,q) \mapsto (c, d)$.
		
		By \threfpart{thm:plys}{4-r} the surface $W^+(12, 1)$ (which by \thref{rmk:4isom} is isomorphic to $W^+(12, 7)$) is birational to the surface given by 
		\begin{equation*}
			\gamma_4^4 - 6 J J' \beta_4 \gamma_4^2 - 16 (J J')^2 \gamma_4 + 3 (J J')^2 (4 J J' -  \beta_4^2) = 0
		\end{equation*}
		in $\mathbb{A}^3$. After making the birational transformation
		\begin{gather*}
			\gamma_4 =  \frac{(2 p^{3} q + 12 p^{2} q^{2} + 18 p q^{3} + 9 q^{4} + 24 p q^{2} + 24 q^{3} - 24 p q - 8 q^{2} - 32 q + 16)^{2}}{64 (p + 4 q - 4) (p + q + 2)^{3}}	\gamma_4' \\
			q = \frac{\left((\gamma_4')^{2} - 3 p + 4 \gamma_4' - 8\right) q' - 2 p \gamma_4' - 4 p + 4 \gamma_4' - 16}{6 \gamma_4'}
		\end{gather*}
		we find that $W^+(12, 1)$ is birational to the quadric surface 
		\begin{equation*}
			3(q')^2 -2 \gamma_4' q' + 8 q' + p + 4 = 0
		\end{equation*}
		in $\bbA^3$. We parametrise this quadric surface by setting
		\begin{equation*}
			(p, q', \gamma_4') = \left( \frac{-4(s - 3)(s + 1)(3s - 2t - 1)^2}{3 (s - 1)^2 (s - 2t - 1)(3s + 2t + 1)}, \frac{8s(t - 1)}{3 (s - 1) (s - 2 t - 1)}, \frac{4 s (s - 3) (3 s - 2 t - 1)}{(s - 1)(s - 2t - 1)(3s + 2t + 1)} \right).
		\end{equation*}
		The forgetful map $W^+(12, 1) \to W(6,1)$ is given by $(s,t) \mapsto (p,q)$.
		
		\subsubsection{The Surface $Z(12, 1)$}
		The $\alpha_4$, $\beta_4$, and $\gamma_4$ are elements of $\bbQ( W^+(12, 1) )$ which satisfy the conditions of \threfpart{thm:plys}{4-r}. Let $\delta_4 = 3(JJ' - (\alpha_4 + 1)^2)$ be as in \threfpart{thm:plys}{4-r}. By construction there exist $\alpha_3, \beta_3, \delta_3 \in \bbQ( W^+(12, 1) )$ which satisfy \threfpart{thm:plys}{3-1}.

		As in Section~\ref{sec:fibprod} let $\mathcal{E}/\bbQ( Z^+(12, 1) )$ be an elliptic curve with $j$-invariant $1728J$. Let $\mathcal{E}_1'/\bbQ( Z^+(12, 1) )$ and $\mathcal{E}_2'/\bbQ( Z^+(12, 1) )$ be elliptic curves with $j$-invariant $1728J'$ which are $(3, 1)$ and $(4,1)$-congruent to $\mathcal{E}$ respectively. 

		It follows from \threfpart{thm:plys}{4-r} and \ref{thm:plys_3-1} that up to square factors $c_6(\mathcal{E})$ is equal to both $3\delta_3c_6(\mathcal{E}_1')$ and $3\alpha_4\delta_4 c_6(\mathcal{E}_2')$. In particular $\mathcal{E}_1'$ is isomorphic to the quadratic twist of $\mathcal{E}_2'$ by $\alpha_4 \delta_4 \delta_3$. 
		
		Therefore by the discussion in Section~\ref{sec:fibprod} the surface $W(12, 1)$ is birational to the double cover of $\bbA^2$ given by
		\begin{equation*}
			w^2 = \alpha_4 \delta_4 \delta_3.
		\end{equation*}
		
		The rational function $\alpha_4 \delta_4 \delta_3$ is equal to $3 s(t - 1) (-s + 2 t + 1) (3 s + 2 t + 1) (s t + 2 s - 2 t - 1)$ up to a square factor. Therefore the surface $W(12, 1)$ is birational to the affine surface given by
    \begin{equation*}
      (w')^2 = 3 s(t - 1) (-s + 2 t + 1) (3 s + 2 t + 1) (s t + 2 s - 2 t - 1).
    \end{equation*}
    We parametrise this surface by setting
		\begin{align*}
			(s,t,w') = \bigg( \frac{-2 u^2 - u v^2 - u + 4 v^2}{u v^2 + u + 4 v^2}, &\frac{2 u^2 v^2 - 4 u^2 + u v^4 - 6 u v^2 - 3 u - 8 v^4}{(2u + v^2 + 3)(u v^2 + u + 4 v^2)}, \\ &\qquad\qquad \frac{36 v (u + 2) (u + v^2 + 1) (u + 2 v^2) (2 u^2 + u v^2 + u - 4 v^2)^2}{(2 u + v^2 + 3)^2 (u v^2 + u + 4 v^2)^3} \bigg).
		\end{align*}
		Hence $W(12, 1)$ is birational to $\bbA^2$ and the forgetful map $W(12, 1) \to W^+(12, 1)$ is given by $(u,v) \mapsto (s,t)$.
		
		Let $F_{12, 1}(u,v)$ be the polynomial from the statement of \thref{thm:Zequations}. It follows immediately that $Z(12, 1)$ is birational to the affine surface given by the equation 
		\begin{equation*}
			z^2 = F_{12, 1}(u,v)
		\end{equation*}
		since $(J - J')^2$ is equal to $F_{12, 1}(u,v)$ up to a square factor. 
		
		\subsubsection{The Surface $Z(12, 7)$}
		By \thref{thm:plys} the surface $W(12, 7)$ is birational to the double cover of $W^+(12, 1)$ given by 
		\begin{equation*}
			w^2 = \alpha_4 \delta_3.
		\end{equation*}
		The rational function $\alpha_4 \delta_3$ is equal to $3 (t - 1) (-9 s^{3} + 12 s^{2} t + 4 s t^{2} + 15 s^{2} + 4 s t - 12 t^{2} + s - 12 t - 3)$ up to a square factor. 
		
		Therefore the surface $W(12, 7)$ is birational to the affine surface given by $(w')^2 = 3 (t - 1) (-9 s^{3} + 12 s^{2} t + 4 s t^{2} + 15 s^{2} + 4 s t - 12 t^{2} + s - 12 t - 3)$. We parametrise this surface by setting
		\begin{equation*}
			(s, t, w') = \left( \frac{(u - 1) (u^{2} - v^{2} + 4 u + 1)}{u^3 + u^2 - u v^2 - 3 u - v^2 + 1}, \frac{u^{3} - u v^{2} + 7 u^{2} - v^{2} + 9 u + 1}{u^3 + u^2 - u v^2 - 3 u - v^2 + 1}, \frac{36 u v (u + 2) (u^2 + 4 u - v^2 + 1)}{(u^3 + u^2 - u v^2 - 3 u - v^2 + 1)^2} \right).
		\end{equation*}
		Hence $W(12, 7)$ is birational to $\bbA^2$ and the forgetful map $W(12, 7) \to W^+(12, 1)$ is given by $(u,v) \mapsto (s,t)$.
		
		Let $F_{12, 7}(u,v)$ be the polynomial from the statement of \thref{thm:Zequations}. It follows immediately that $Z(12, 7)$ is birational to the affine surface given by the equation 
		\begin{equation*}
			z^2 = F_{12, 7}(u,v) 
		\end{equation*}
		since $(J - J')^2$ is equal to $F_{12, 7}(u,v)$ up to a square factor. 
	
	\subsection{Models for \texorpdfstring{$Z(12, 5)$}{Z(12,5)} and \texorpdfstring{$Z(12, 11)$}{Z(12,11)}}	\label{sec:Z12_5_11}
	
		For each $r = 1, 3$ let $\mathcal{W}(4, r)$ be the surface in $\bbA^5$ defined by the polynomials in \threfpart{thm:plys}{2-1}-\ref{thm:plys_4-r}. The rational map $\bbA^2 \dashrightarrow \mathcal{W}(4, r)$ given by taking 
		\begin{gather*}
			JJ' = \frac{16(-2 a b + b^{2} + 2 a + 2 b)^{3}}{27(b - 1)^{2}(-a + b + 1)^{4}},	\\
			(J - 1)(J' - 1) = \frac{(a^{2} b + 2 a b^{2} + b^{3} - a^{2} + b^{2} - 2 a + b + 1)^{2}}{(b - 1)^{2} (-a + b + 1)^{4}},
		\end{gather*}
		(and solving for $\alpha$, $\beta$, and $\gamma$) is birational. In particular $W(4, r)$ is birational to $\bbA^2$ and the forgetful maps $W(4, r) \to W(1)$ are given by the rational functions $JJ'$ and $(J-1)(J'-1)$ above.
		
		Let $\Delta = \Delta(\mathcal{E})$ and $\Delta' = \Delta(\mathcal{E}')$ be the discriminants of elliptic curves $\mathcal{E}/\bbQ(Z(4,1))$ and $\mathcal{E}'/\bbQ(Z(4,1))$ with $j$-invariants $1728J$ and $1728J'$ respectively. We write $W(4, 1)^{\sqrt[3]{\Delta\Delta'}}$ for the triple cover of $W(4, 1)$ whose function field is  given by adjoining $\sqrt[3]{\Delta \Delta'}$ to $\bbQ(W(4,1))$.
		
		Then $W(4, 1)^{\sqrt[3]{\Delta\Delta'}}$ is birational to the surface given by
		\begin{equation*}
			\alpha_3^3 - JJ' = 0
		\end{equation*}
		in $\bbA^3$. By making a change of coordinates $\alpha_3 = \frac{2(-2 a b + b^{2} + 2 a + 2 b)}{3(b - 1)(-a + b + 1)^{2}} \alpha_3'$ we see that $W(4, 1)^{\sqrt[3]{\Delta\Delta'}}$ is birational to the surface
		\begin{equation*}
			(\alpha_3')^3 - 2(b - 1)(-a + b + 1)^2 = 0
		\end{equation*}
		in $\bbA^3$. We parametrise this cubic surface by setting
		\begin{equation*}
			(a,b,\alpha_3') = \left(\frac{-c d^3 + 8 d^3 + 4}{2 (c + 2d^3 + 1)}, \frac{- c + 2 d^3 + 1}{c + 2d^3 + 1}, \frac{-cd^{2}}{c + 2 d^3 + 1} \right).
		\end{equation*}
		Hence $W(4, 1)^{\sqrt[3]{\Delta\Delta'}}$ is birational to $\bbA^2$ and the forgetful map $W(4, 1)^{\sqrt[3]{\Delta\Delta'}} \to W(4, 1)$ is given by $(c,d) \mapsto (a,b)$.
		
		By \threfpart{thm:plys}{3-2} the surface $W^+(12, 5)$ (which by \thref{rmk:4isom} is isomorphic to $W^+(12, 11)$) is birational to the surface given by 
		\begin{equation*}
			\beta_3^4 -6 (\alpha_3 + 1) ( J - 1) (J' - 1) \beta_3^2 -8 (( J - 1) (J' - 1))^2 \beta_3 - 3 ((\alpha_3 - 1) ( J - 1) (J' - 1))^2 = 0
		\end{equation*}
		in $\mathbb{A}^3$. 
		This is clearly birational to the surface given by
		\begin{gather*}
			\beta_3^2 - \theta = 0 \\
			\theta^2 -6 (\alpha_3 + 1) ( J - 1) (J' - 1) \theta -8 (( J - 1) (J' - 1))^2 \beta_3 - 3 ((\alpha_3 - 1) ( J - 1) (J' - 1))^2 = 0
		\end{gather*}
		in $\bbA^4$. To simplify the equations for this surface, we view the forms above as defining a quadric intersection in $\bbA^2_{\bbQ(c,d)}$. By using `minimisation' techniques similar to those in~\cite[Section~4.3]{CFS_MARO234COEC}, but applied over the base ring $\bbZ[c,d]$, we arrive at the change of coordinates
		\begin{gather*}
			\theta = \frac{g_3(c,d)^2 (g_1(c,d) \theta' + g_2(c,d) \beta_3' + 2(2d^3 + 1)g_2(c,d))}{3(\beta_3' + 2(2d^3 + 1))}, \\ 
			\beta_3 = \frac{g_3(c,d) \left(-3(d^2 c + 4(2d^3 + 1))\beta_3' + 2(2d^3 + 1)(2c^2 -3(d^2 + 4)c +6(2d^3 + 1)\right)}{3(\beta_3' + 2(2d^3 + 1))}
		\end{gather*}
		where
		\begin{align*}
			g_1(c,d) &= 8 (2d^3 + 1)(c^2 - 6c + 9(2d^3 + 1)) ,\\
			g_2(c,d) &= (3d^4 - 8d^3 - 4) c^2 + 24(d^2 + 1)(2d^3 + 1)c + 12(2d^3 + 1)^2 ,\\
			g_3(c,d) &= \frac{d^6 c^3 - 8(2d^3 + 1)^2 c^2 + 48(2d^3 + 1)^2 c - 8(2d^3 + 1)^3}{c^4d^8} .
		\end{align*}
		It follows that $W^+(12, 5)$ is birational to the surface given by the pair of equations
		\begin{gather*}
			2 (\beta_3')^{2} - 2 \beta_3' \theta' - (\theta')^{2} + \left(-8 d^{3} - 4\right) \beta_3' + \left(-8 d^{3} - 4\right) \theta' + c^{2} d^{4} - 8 c d^{5} - 4 c d^{2} = 0, \\
			-3 (\beta_3')^{2} + 3 (\theta')^{2} -6 c d^{2} \beta_3' - 3 c^{2} d^{4} + 8 c^{2} d^{3} - 48 c d^{3} + 4 c^{2} - 24 c = 0
		\end{gather*}
		in $\bbA^4$. We parametrise this surface by setting
		\begin{equation*}
			(c,d,\beta_3',\theta') = \left( \frac{6(s + t) h_1(s,t)}{h_3(s,t)}, \frac{t}{s + t + 2}, \frac{2 h_1(s,t)h_2(s,t) }{(s + t +	2)^2 h_3(s,t)} \frac{-4 (s - 1) (s^2 + s t + s + 2 t + 1) h_1(s,t)}{ (s + t + 2)^2 h_3(s, t)} \right)
		\end{equation*}
		where
		\begin{align*}
			h_1(s,t) &= s^{3} + \left(3 t + 6\right) s^{2} + \left(3 t^{2} + 12 t + 12\right) s + 3 t^{3} + 6 t^{2} + 12 t + 8, \\
			h_2(s,t) &= 4 s^3 + 4 s^2 t - 3 s t^2 - 2 s t - 6 s - 3 t^3 - 2 t + 2, \\
			h_3(s,t) &= 2 s^{4} + \left(4 t + 4\right) s^{3} + \left(5 t^{2} + 18 t + 12\right) s^{2} + \left(6 t^{3} + 20 t^{2} + 24 t + 10\right) s 
					 \\&\quad + 3 t^{4} + 6 t^{3} + 11 t^{2} + 8 t - 1.
		\end{align*}
		Therefore $W^+(12, 5)$ is birational to $\bbA^2$ and the forgetful map $W^+(12, 5) \to W(4, 1)^{\sqrt[3]{\Delta\Delta'}}$ is given by $(s,t) \mapsto (c,d)$.
		
		\subsubsection{The Surface $Z(12, 5)$}

		Note that $\alpha_3, \beta_3 \in \bbQ( W^+(12, 5) )$ satisfy the conditions of \threfpart{thm:plys}{3-2}. By construction there exist elements $\alpha_4$, $\beta_4$, $\gamma_4$ of $\bbQ( W^+(12, 1) )$ which satisfy \threfpart{thm:plys}{2-1}--\ref{thm:plys_4-r}. Let $\delta_4 = 3(JJ' - (\alpha_4 + 1)^2)$ be as in \threfpart{thm:plys}{4-r}
		
		By \thref{thm:plys} the surface $W(12, 5)$ is birational to the double cover of $W^+(12, 5)$ given by 
		\begin{equation*}
			w^2 = \alpha_4 \delta_4 \beta_3.
		\end{equation*}
		The rational function $\alpha_4 \delta_4 \beta_3$ is equal to  $-s^2 + t^2 + 1$ up to a square factor. 
		
		Therefore the surface $W(12, 5)$ is birational to the affine surface $(w')^2 = -s^2 + t^2 + 1$. We parametrise this surface by setting
		\begin{equation*}
			(s, t, w') = \left( \frac{-u^{2} + v^{2} + 2 u - 2}{u^{2} - 2u - v^2}, \frac{2u - 2}{u^2 - 2u - v^2}, \frac{2v}{u^2 - 2 u - v^2} \right).
		\end{equation*}
		Hence $W(12, 5)$ is birational to $\bbA^2$ and the forgetful map $W(12, 5) \to W^+(12, 5)$ is given by $(u,v) \mapsto (s,t)$.
		
		Let $F_{12, 5}(u,v)$ be the polynomial from the statement of \thref{thm:Zequations}. It follows immediately that $Z(12, 5)$ is birational to the affine surface given by the equation 
		\begin{equation*}
			z^2 = F_{12, 5}(u,v) 
		\end{equation*}
		since $(J - J')^2$ is equal to $F_{12, 5}(u,v)$ up to a square factor. 
		
		\subsubsection{The Surface $Z(12, 11)$}
		By \thref{thm:plys} the surface $W(12, 11)$ is birational to the double cover of $W^+(12, 5)$ given by 
		\begin{equation*}
			w^2 = \alpha_4 \beta_3.
		\end{equation*}
		The rational function $\alpha_4 \beta_3$ is equal to $ -(s + t)(s + t + 2)$ up to a square factor.
		
		Therefore the surface $W(12, 11)$ is birational to the affine surface $(w')^2 = -(s + t)(s + t + 2)$. We parametrise this surface by
		\begin{equation*}
			(s,t,w') = \left(\frac{-u v^2 + u - 3 v^2 - 1}{u v^2 + u + v^2 + 1}, \frac{-u + 1}{u + 1}, \frac{2 v}{v^2 + 1} \right).
		\end{equation*}
		Hence $W(12, 11)$ is birational to $\bbA^2$ and the forgetful map $W(12, 11) \to W^+(12, 5)$ is given by $(u,v) \mapsto (s,t)$.
		
		Let $F_{12, 11}(u,v)$ be the polynomial from the statement of \thref{thm:Zequations}. It follows immediately that $Z(12, 11)$ is birational to the affine surface given by the equation 
		\begin{equation*}
			z^2 = F_{12, 11}(u,v) 
		\end{equation*}
		since $(J - J')^2$ is equal to $F_{12, 11}(u,v)$ up to a square factor.  

\section{Curves on \texorpdfstring{$Z(12,r)$}{Z(12,r)} and examples of \texorpdfstring{$12$}{12}-congruences} \label{sec:curves}

	We now identify some curves of small genus on  our models for $Z(12, r)$ and we prove \thref{thm:infmany}. 
	
	\subsection{The modular curves \texorpdfstring{$X_0(m)$}{X0(m)} on \texorpdfstring{$Z(12,r)$}{Z(12,r)}} \label{sec:X0m_on_Z}
		Let $m$ be an integer coprime to $N$, and suppose that $m$ is equal to $r$ modulo $N$ (up to a square factor). The (non-compact) modular curve $Y_0(m)$ may be naturally embedded in $Z(N,r)$. More precisely we may regard a non-cuspidal $K$-point on $Y_0(m)$ as a triple $(E, E', \phi)$ where $E/K$ and $E'/K$ are elliptic curves and $\phi$ is a cyclic $m$-isogeny. In particular we have a natural map 
		\begin{align*}
			Y_0(m) & \stackrel{\psi}{\longrightarrow} Z(N, m)\\
			(E, E', \phi) &\mapsto (E, E', \phi \vert_{E[N]} )
		\end{align*}
		defined over $\bbQ$. A pair of points $(E_1, E_1', \phi_1)$ and $(E_2, E_2', \phi_2)$ are identified under this morphism if and only if the the pairs $(E_1, E_1')$ and $(E_2, E_2')$ are simultaneous quadratic twists, in which case they define the same point on $Y_0(m)$. In particular $\psi$ is an embedding.
		
		\begin{remark}
			In fact, let $P = (E, E', \phi)$ be a $K$-point on $Z(N, r)$. Then for each integer $k$ such that $k^2 \equiv 1 \pmod{N}$ we obtain another $K$-point $P_a = (E, E', \phi \circ [k])$ on $Z(N, r)$ which is equal to $P$ if and only if $k$ is equal to either $1$ or $-1$ modulo $N$. 
			
			In particular when $N = 12$ points on $Z(N, r)$ come in pairs, namely $P_1$ and $P_5$ in the above notation, which have the same image under $Z(12, r) \to Z^+(12, r)$. Hence each modular curve $Y_0(m)$ as above appears twice on $Z(12, m)$. In the notation of Kani--Schanz~\cite{KS_MDQS} these are the ``Hecke curves'' $T_{m, 1}$ and $T_{m,5}$.  
		\end{remark}
		
		The image of $Y_0(m)$ under the quotient $Z(N, r) \to W(N, r)$ is $Y_0^+(m) \colonequals Y_0(m)/w_m$, where $w_m$ denotes the Fricke involution on $X_0(m)$.
		
		We identify the curves $X_0(m)$ on our models for $Z(N, r)$ for each $m$ such that $X_0(m)$ appears in \texttt{Magma}'s \texttt{SmallModularCurve} database. The equations for these curves are recorded in Tables~\ref{table:X0_Z121}--\ref{table:X0_Z1211}, the two choices of sign give the pair of curves $T_{m,1}$ and $T_{m, 5}$. 
		
		Note that the modular curves $X_0(5)$ on $Z(12, 5)$ and $X_0(7)$ on $Z(12, 7)$ have been blown down (this is the case since we determined $Z(12, r)$ only up to birational equivalence). Equations for these curves on a blowup of our models are given by 
		\begin{equation*}
			z^2 = F_{12,5} \left(\tfrac{1}{2}(3 - \varepsilon),  \tfrac{1}{2} (\pm 1 \pm \varepsilon + \varepsilon^2 t) \right) = \tfrac{1}{4} \left(-t^2 \mp 44 t + 16 \right) \varepsilon^4 + O(\varepsilon^5)
		\end{equation*}
		and 
		\begin{equation*}
			z^2 = F_{12,7}(-\varepsilon , \pm1 \pm \varepsilon + \varepsilon^2 t) = (4t^2 \pm 52t - 27)\varepsilon^4 + O(\varepsilon^5)
		\end{equation*}
		respectively.

		\begin{table}[H]
			\centering 
			\begin{tabular}{c|c}
				$m$			& Equation for $X_0^+(m)$ on $W(12, 1)$ \\
				\hline
				$13$		& $\pm v + 1$	\\
				$25$		& $2 u + v^2 \pm 2 v + 1$	\\
				$37$		& $2 u^2 + (v^2 \pm 2 v + 1)u \pm v^3 - 3 v^2 \pm 3 v - 1$	\\
				$49$		& $4u^2 + 2(3v^2 \pm 2v + 3)u + v^4 \pm 4v^3 + 6v^2 \pm 4v + 1$	\\
				$61$		& $4 u^3 + 2( v^2 \pm 2 v + 5)u^2 + (-v^4 \pm 6 v^3 - 4 v^2 \pm 10 v + 5) u \pm v^5 - 5v^4 \pm 10v^3 - 10v^2 \pm 5v - 1$
					\\
			\end{tabular}
      \caption{Modular curves on $Z(12, 1)$.} \label{table:X0_Z121}
		\end{table}
	
		\begin{table}[H]
			\centering 
			\begin{tabular}{c|c}
				$m$			& Equation for $X_0^+(m)$ on $W(12, 5)$ \\
				\hline
				$17$		&	$u \pm v$ \\
				$29$		&	$u^2 + 2 u - v^2 \pm 2v - 4$\\
				$41$		&	$u^3 - (\pm v - 2) u^2 + (-v^2 - 12) u \pm v^3 - 2 v^2 \pm 4 v + 8$\\
				$53$		&	$u^4 - (2 v^2 \pm 4 v) u^2 + (\pm 8 v - 16) u + v^4 \pm 4 v^3 + 8 v^2 + 16$\\
			\end{tabular}
      \caption{Modular curves on $Z(12, 5)$.} \label{table:X0_Z125}
		\end{table}
	
		\begin{table}[H]
			\centering 
			\begin{tabular}{c|c}
				$m$			& Equation for $X_0^+(m)$ on $W(12, 7)$ \\
				\hline
				$19$		&	$1 \pm v$ \\
				$31$		&	$u^2 - (\pm v - 1) u \pm 2 v - 2$ \\
				$43$		&	$-(\pm v - 1) u^2 + (v^2 \pm 2 v + 1) u + 2 v^2 - 2$\\
				$55$		&	$(v^2 \mp 6 v + 1) u^2 + (-2 v^2 \pm 4 v - 2) u - v^4 \mp 2 v^3 \pm 2 v + 1$\\
			\end{tabular}
      \caption{Modular curves on $Z(12, 7)$.} \label{table:X0_Z127}
		\end{table}
	
		\begin{table}[H]
			\centering 
			\begin{tabular}{c|c}
				$m$			& Equation for $X_0^+(m)$ on $W(12, 11)$ \\
				\hline
				$11$		&	$1 \pm v$\\
				$23$		&	$\pm v u + v^2 \pm v + 1$ \\
				$35$		&	$-(\pm v^3 \pm 3 v) u - v^4 \pm 3 v^3 \pm v + 1$\\
				$47$		&	$(\pm 3 v^3 + 2 v^2 \pm 3 v) u + v^4 \pm v^3 + 4 v^2 \pm v + 1$\\
				$59$		&	$(-3 v^4 \pm 3 v^3 - v^2 \pm v) u \pm v^5 + 8 v^4 \pm 7 v^3 + 11 v^2 \pm 4 v + 1 $\\
				$71$		&	$(\pm 7 v^5 \pm 10 v^3 \mp v) u + v^6 \mp v^5 + 15 v^4 \pm 10 v^3 + 15 v^2 \pm 7 v + 1$ 
			\end{tabular}
      \caption{Modular curves on $Z(12, 11)$.} \label{table:X0_Z1211}
		\end{table}

	\subsection{Infinite families of \texorpdfstring{$12$}{12}-congruent elliptic curves and examples over \texorpdfstring{$\bbQ(t)$}{Q(t)}} \label{sec:inffams}		
		We now prove \thref{thm:infmany}. Specifically, we identify curves $C \subset Z(12, r)$ of genus $0$ or $1$ with infinitely many rational points. We use these to construct pairs of elliptic curves defined over $\bbQ(C)$ which are $(12, r)$-congruent. Once we have such a pair of elliptic curves the following lemma will allow us to deduce \thref{thm:infmany}. 
		
		Recall that if $E/\bbQ$ is an elliptic curve, then fixing a basis for $E[N]$ we obtain an isomorphism $\Aut(E[N]) \cong \GL_2(\bbZ/N\bbZ)$ and hence a representation 
		\begin{equation*}
			\bar{\rho}_{E, N} \colon G_\bbQ \to \GL_2(\bbZ/N\bbZ)
		\end{equation*}
		which we call the \emph{mod $N$ Galois representation attached to $E/\bbQ$} (for more details see \cite[Section~2.2]{RSZB_LAIOGFECOQ}).
	 
		\begin{lemma} \thlabel{lemma:galoisrep}
			Let $C/\bbQ$ be a smooth projective curve with infinitely many rational points and let $\mathcal{E}/\bbQ(C)$ and $\mathcal{E}'/\bbQ(C)$ be elliptic curves which are not both isotrivial. For a point $P \in C(\bbQ)$ let $\mathcal{E}_P/\bbQ$ be the specialisation of $\mathcal{E}$ at $P$ (and similarly for $\mathcal{E}'$). 
			\begin{enumerate}[label=\eniii]
				\item Suppose there exists a point $P \in C(\bbQ)$ such that $\mathcal{E}_P$ and $\mathcal{E}'_P$ are nonsingular and do not have CM. If $j(\mathcal{E}_P)$ and $j(\mathcal{E}'_P)$ are not $j$-invariants of $m$-isogenous elliptic curves over $\bbQ$ for any integer $m$ then for all but finitely many $Q \in C(\bbQ)$ we have that $j(\mathcal{E}_Q)$ and $j(\mathcal{E}'_Q)$ are not the $j$-invariants of $m$-isogenous elliptic curves for any integer $m$.

				\item Suppose there exists a point $P \in C(\bbQ)$ such that $\mathcal{E}_P$ is nonsingular and the mod $N$ Galois representation attached to $\mathcal{E}_P/\bbQ$ is surjective. Then there exist infinitely many $Q \in C(\bbQ)$ such that the mod $N$ Galois representation attached to $\mathcal{E}_Q/\bbQ$ is surjective.
			\end{enumerate}
			
			\begin{proof}
				(i) This follows from a similar argument to \cite[Theorem~1.5]{F_OFO9CEC} which we reproduce here for completeness. 
				
				First suppose $E$ and $E'$ are non-CM elliptic curves defined over $\bbQ$ and $\phi \colon E \to E'$ is an isogeny defined over $\Qbar$. Then $\widehat{\phi} \phi = [n] \in \End(E)$ for some $n \in \bbZ$ and by comparing degrees for any $\sigma \in G_{\bbQ}$ we have $(\widehat{\phi})^\sigma \phi = \chi(\sigma) [n] \in \End(E)$ for some character $\chi \colon G_\bbQ \to \{\pm 1\}$. Let $d \in \bbQ$ be such that $\bbQ(\sqrt{d})$ is cut out by the kernel of $\chi$. Let $E^d/\bbQ$ be the quadratic twist of $E$ by $d$ and let $t \colon E^d \to E$ be a $\Qbar$-isomorphism. Define $\phi' = \phi \circ t$, then for each $\sigma \in G_\bbQ$ we have $(\widehat{\phi}')^\sigma \phi' = [n] = \widehat{\phi}' \phi'$. But $\phi'$ is a surjection on $\Qbar$-points so $(\widehat{\phi}')^\sigma = \widehat{\phi}'$ for all $\sigma \in G_{\bbQ}$. Therefore $\phi' \colon E^d \to E'$ is a $\bbQ$-isogeny. That is, $j(E)$ and $j(E')$ are the $j$-invariants of $\bbQ$-isogenous elliptic curves. 
				
				It therefore suffices to show that there are at most finitely many $Q \in C(\bbQ)$ such that $j(\mathcal{E}_Q)$ and $j(\mathcal{E}'_Q)$ are $j$-invariants of $\bbQ$-isogenous elliptic curves.
				
				Let $C'$ be the image of the map $C \to X(1) \times X(1)$ given by extending $t \mapsto (j(\mathcal{E}), j(\mathcal{E}'))$. Consider the curves, $T_m$, on $X(1) \times X(1)$ which are the image of the modular curves $X_0(m)$ under the graph of the Hecke correspondence. Since the curves $T_m$ are geometrically irreducible and $j(\mathcal{E}_P)$ and $j(\mathcal{E}'_P)$ are not $j$-invariants of $\bbQ$-isogenous elliptic curves, the curve $C'$ meets each $T_m$ at finitely many points (over $\Qbar$). 
				
				Mazur~\cite{M_RIOPD} and Kenku~\cite{K_OTNOQICOECIEQIC} showed that if an elliptic curve admits a cyclic $m$-isogeny over $\bbQ$ then $m \leq 163$. In particular there are finitely many rational points on $C'$ which lie on the union of all $T_m$. The claim follows since the map $j \colon C \to C'$ has finite degree.
				
				(ii) We define a morphism $j \colon C \to \bbP^1 \cong X(1)$ by extending $j(\mathcal{E})$. Let $\mathcal{M}$ be the set of all maximal subgroups of $\GL_2(\bbZ/N\bbZ)$. For each $H \in \mathcal{M}$ let $X(H)/\bbQ$ denote the modular curve whose $K$-points parametrise elliptic curves $E/K$ such that $\bar{\rho}_{E,N}(G_\bbQ)$ is contained in a subgroup of $\GL_2(\bbZ/N\bbZ)$ conjugate to $H$ (see e.g., \cite[Section~2]{RSZB_LAIOGFECOQ}). Let $C_H$ be the fibre product 
				\[\begin{tikzcd}
					{C_H} & {X(H)} \\
					{C} & {X(1)}
					\arrow["\phi_H"', from=1-1, to=2-1]
					\arrow["j", from=2-1, to=2-2]
					\arrow[from=1-1, to=1-2]
					\arrow[from=1-2, to=2-2]
				\end{tikzcd}\]
				We claim that the set 
				\begin{equation*}
					C(\bbQ) \setminus \bigcup_{H \in \mathcal{M}} \phi_H(C_H(\bbQ))
				\end{equation*}
				is infinite. 
				
				First suppose that $C$ has genus $0$. Note that the morphisms $\phi_H$ have degree greater than $> 1$ because $\mathcal{E}_P/\bbQ$ has surjective mod $N$ Galois representation.  Therefore since $\mathcal{M}$ is finite $\bigcup_{H \in \mathcal{M}} \phi_H(C_H(\bbQ))$ is a thin set in the sense of Serre. The claim follows from Hilbert's Irreducibility Theorem \cite[Proposition~3.4.2]{S_TIGT}.
				
				If $C$ has genus $1$ we make $C$ into an elliptic curve by declaring $P$ to be the identity. By Faltings' theorem any of the curves $C_H$ of genus $\geq 2$ have finitely many rational points. Therefore we may assume without loss of generality that for each $H$ the curve $C_H$ has genus $1$ (and that $C_H(\bbQ)$ is non-empty). 
				
				Since every morphism between elliptic curves is the composition of a translation and an isogeny \cite[Example~4.7]{S_TAOEC} the images $\phi_H(C_H(\bbQ))$ are cosets of subgroups $\Phi_H \subset C_H(\bbQ)$ which have finite index by the weak Mordell--Weil theorem. Because $\mathcal{E}_P/\bbQ$ has surjective mod $N$ Galois representation the point $P$ is not contained in $\phi_H(C_H(\bbQ))$ for any $H$. 
				
				Since we took $P$ to be the identity, for each $H$ the image $\phi_H(C_H(\bbQ))$ is a non-trivial coset of $\Phi_H$. In particular the subgroup $\bigcap_H \Phi_H \subset C(\bbQ)$ is disjoint from $\phi_H(C_H(\bbQ))$ for every $H$. The claim follows by noting that the intersection $\bigcap_H \Phi_H$ has finite index in $C(\bbQ)$, and is therefore infinite.
			\end{proof}
		\end{lemma}

		\begin{example}[A family of $(12,1)$-Congruences] \thlabel{ex:121}
			Examples of pairs of $(12, 1)$-congruent elliptic curves over $\bbQ(t)$ have previously been given by Chen~\cite[Section~7.3]{C_COEC} and Fisher~\cite{F_EMSFCOEC}. We include one here for completeness. The model $z^2 = F_{12, 1}(u,v)$ has an obvious elliptic fibration over the $t$-line given by $y^2 = F_{12, 1}(x, t)$. This model has a pair of sections with $x$-coordinates $0$. By putting this fibration in Weierstrass form we see that $Z(12, 1)$ is birational to the elliptic K3 surface with Weierstrass equation
			\begin{align*}
				y^2 &= x^3 - 27(t^8 + 46 t^6 + 859 t^4 - 186 t^2 + 9)x \\&\quad- 54 (t^{12} + 69 t^{10} + 2082 t^8 + 24731 t^6 - 7848 t^4 + 621 t^2 + 
				27). 
			\end{align*}
			Moreover the sections with $x$-coordinate $\frac{3(2 t^8 + 37 t^6 + 450 t^4 - 27 t^2 + 54)}{(t^2 - 3)^2}$ have infinite order. Taking multiples of this section gives an infinite family of examples over $\bbQ(t)$ (a different family is given by Fisher~\cite[Corollary~1.3]{F_EMSFCOEC}).
		
			In order to give a family where the $12$-congruent elliptic curves have simple Weierstrass equations we instead note that the curve on $Z(12, 1)$ given by the vanishing of $2u + v^2 + 3v$ has genus $0$. By parametrising this curve we obtain a pair of $(12, 1)$-congruent elliptic curves $\mathcal{E}_{12, 1}/\bbQ(t)$ and $\mathcal{E}_{12, 1}'/\bbQ(t)$ where
			\begin{align*}
				\mathcal{E}_{12, 1}  &: y^2 = x^3 - 3(t^2 + 1)(4t^2 + 2t + 1)p(t)x - 2(t^2 + 1)^2(4t^2 + 2t + 1)q(t)\\
				\mathcal{E}_{12, 1}' &: y^2 = x^3 - 3(t^2 + 1)(4t^2 + 2t + 1)p'(t)x + 2(t^2 + 1)^2(4t^2 + 2t + 1)q'(t)
			\end{align*}
			and
			\begin{flalign*}
				p(t) &= 36 t^{12} + 234 t^{11} + 693 t^{10} + 1530 t^9 + 2808 t^8 + 4128 t^7 + 4964 t^6 + 4792 t^5 + 3652 t^4 + 2000 t^3 \\&\quad + 736 t^2 + 160 t + 16, \\
				q(t) &= 432 t^{18} + 4320 t^{17} + 20196 t^{16} + 63720 t^{15} + 159111 t^{14} + 330453 t^{13} +581103 t^{12} + 875862 t^{11} \\&\quad + 1137762 t^{10} + 1270440 t^9 + 1208628 t^8 +960528 t^7 + 622592 t^6 + 320528 t^5 + 127472 t^4 \\&\quad + 37728 t^3 + 7840 t^2 +1024 t + 64, \\
				p'(t) &= 55396 t^{12} - 97238 t^{11} + 230581 t^{10} - 206878 t^9 + 177280 t^8 - 86384 t^7 + 2612 t^6 - 6600 t^5 \\&\quad + 900 t^4 - 1008 t^3 + 288 t^2 + 288 t + 144	, \\
				q'(t) &= 26076112 t^{18} - 62137376 t^{17} + 165178492 t^{16} - 207914808 t^{15} + 236336753 t^{14}- 163049017 t^{13} \\&\quad + 71325401 t^{12} - 36145662 t^{11} + 7345854 t^{10} - 7544088 t^9 + 2972844 t^8 + 1350576 t^7 + 531072 t^6 \\&\quad - 137808 t^5 - 17712 t^4 - 4320 t^3 + 14688 t^2 + 6912 t + 1728.
			\end{flalign*}
		\end{example}

		\begin{example}[A family of $(12,5)$-Congruences] \thlabel{ex:125}
			Our model for the surface $Z(12, 5)$ contains the genus $1$ curve given by the vanishing of $u^3 + 3u^2 - (v^2 + 10) u - 3 v^2 + 6$. This curve is birational over $\bbQ$ to the elliptic curve $\mathcal{C}/\bbQ$ given by the Weierstrass equation
			\begin{equation*}
				\mathcal{C} : \eta^2 = \xi^3 - \xi^2 - 16\xi + 16.
			\end{equation*}
			The elliptic curve $\mathcal{C}$ has LMFDB label \LMFDBLabel{240.a3} and rank $1$.

			Using the moduli interpretation for $Z(12, 5)$ it is possible to give an explicit pair of $(12, 5)$-congruent elliptic curves over $\bbQ(\mathcal{C})$ in terms of $\xi$ and $\eta$. The Weierstrass equations are too complicated to reproduce here, but we include them in the electronic data~\cite{ME_ELECTRONIC_Z12}.
		
		\end{example}

		\begin{example}[A family of $(12,7)$-Congruences] \thlabel{ex:127}
			Our model for the surface $Z(12, 7)$ contains the genus $0$ curve given by the vanishing of $u^2 + 3 u - 2 v^2 + 2$. By parametrising this curve we obtain a pair of $(12, 7)$-congruent elliptic curves $\mathcal{E}_{12, 7}/\bbQ(t)$ and $\mathcal{E}_{12, 7}'/\bbQ(t)$ where
			\begin{align*}
				\mathcal{E}_{12, 7}  &: y^2 = x^3 - 3 \frac{(t^{2} + 3) (t^{4} - 3 t^2 + 9) p(t)}{(t^{4} + 3)(t^{4} + 27)} x - 2 \frac{(t^{2} + 3) (t^{4} - 3 t^2 + 9) q(t)}{(t^{4} + 3)(t^{4} + 27)} 
				\
			\end{align*}
			\begin{align*}
				\mathcal{E}_{12, 7}' :y^2 = x^3 -3 (t^{2} + 3) (t^{4} - 3 t^2 + 9) &(t^{4} + 3) (t^{4} + 27) p'(t) x \\&-2  (t^{2} + 3) (t^{4} - 3 t^2 + 9) (t^{4} + 3)^{2} (t^{4} + 27)^{2} q'(t)
			\end{align*}
			and
			\begin{flalign*}
				p(t) &= t^{18} - 18 t^{16} + 135 t^{14} - 729 t^{12} + 2889 t^{10} - 7533 t^{8} + 15093 t^{6} - 15795 t^{4} + 53946 t^{2} + 2187, \\
				q(t) &= t^{26} - 27 t^{24} + 309 t^{22} - 2133 t^{20} + 10341 t^{18} - 35559 t^{16} + 85158 t^{14} - 109350 t^{12} - 141345 t^{10} \\ &\quad  +1337715 t^{8} - 3133971 t^{6} + 4183731 t^{4} + 3483891 t^{2} - 59049, \\
				p'(t) &= t^{18} + 222 t^{16} - 585 t^{14} + 5031 t^{12} - 22599 t^{10} + 78003 t^{8} - 177147 t^{6} + 295245 t^{4} - 354294 t^{2} + 177147, \\
				q'(t) &= t^{26} - 531 t^{24} - 5739 t^{22} + 38691 t^{20} - 148635 t^{18} + 141345 t^{16} + 984150 t^{14} - 6897798 t^{12} + 25922511 t^{10} \\ &\quad -67847301 t^{8} + 125951517 t^{6} - 164215269 t^{4} + 129140163 t^{2} - 43046721.
			\end{flalign*}
		\end{example}

		\begin{example}[A family of $(12,11)$-Congruences] \thlabel{ex:1211}
			Our model for the surface $Z(12, 11)$ contains a pair of genus $0$ curves given by the vanishing of $(9 v^4 + 30 v^2 + 5) u - 45 v^4 - 6 v^2 + 7$. By parametrising one of these curves we obtain a pair of $(12, 11)$-congruent elliptic curves $\mathcal{E}_{12, 11}/\bbQ(t)$ and $\mathcal{E}_{12, 11}'/\bbQ(t)$ where
			\begin{flalign*}
				\mathcal{E}_{12, 11} : y^2 = x^3 -3& \frac{(25 t^6 + 63 t^4 + 27 t^2 - 27) (t^6 + 15 t^4 + 3t^2 - 27)}{(t^4 - 18t^2 - 27)(t^4 + 6t^2 - 3)} x \\& \qquad \qquad \qquad- 2 \frac{(25 t^6 + 63 t^4 + 27 t^2 - 27)(5t^8 -24t^6 -78t^4 + 81)}{(t^4 - 18t^2 - 27)(t^4 + 6t^2 - 3)}
			\end{flalign*}
			\begin{flalign*}
				\mathcal{E}_{12, 11}'  &: y^2 = x^3 -3 (t^4 - 18 t^2 - 27) (t^4 + 6 t^2 - 3) p'(t) x -2 \frac{(t^4 - 18 t^2 - 27)^2 (t^4 + 6 t^2 - 3)^2 q'(t)}{(25 t^6 + 63 t^4 + 27 t^2 - 27)}			
			\end{flalign*}
			and 
			\begin{flalign*}
				p'(t) &= 25 t^{20} + 18 t^{18} - 1971 t^{16} - 16200 t^{14} - 138942 t^{12} - 367092 t^{10} -
				500526 t^8 -332424 t^6 - 13851 t^4 \\&\quad + 1458 t^2 + 6561, \\
				q'(t) &= 3125 t^{32} + 30000 t^{30} + 94800 t^{28} + 878688 t^{26} + 25070580 t^{24} +
				173317968 t^{22} + 693690912 t^{20} \\&\quad + 1821351744 t^{18} + 3069524646 t^{16} +
				3140007120 t^{14} + 1733188752 t^{12} + 288404064 t^{10} - 233440380 t^8 \\&\quad -
				163447632 t^6 + 7558272 t^4 - 1594323.
			\end{flalign*}
		\end{example}

		We now complete the proof of \thref{thm:infmany}.

		\begin{proof}[Proof of \thref{thm:infmany}]
			Consider the pairs of elliptic curves, $E_{12,r}/\bbQ$ and $E_{12, r}'/\bbQ$ given by specialising the pairs in Examples~\ref{ex:121}, \ref{ex:127}, and \ref{ex:1211} at $t = 2$ and the pair in \thref{ex:125} at the point $(\xi, \eta) = (0, 4)$. The hypotheses of condition (i) in \thref{lemma:galoisrep} are satisfied since for each $r$ the curve $E_{12, r}$ is alone in its $\bbQ$-isogeny class. The first claim then follows from \thref{lemma:galoisrep}(i).

			Note that if the mod $N$ Galois representation attached to an elliptic curve $E/\bbQ$ is not surjective then there exists a maximal subgroup $H \subset \GL_2(\bbZ/N\bbZ)$ such that $\bar{\rho}_{E, N}(G_\bbQ)$ is contained in a subgroup conjugate to $H$. In particular if $p$ is a prime number not dividing $N$ or the denominator of $j(E)$ then the reduction of $j(E)$ modulo $p$ is the $j$-invariant of an $\bbF_p$-rational point on $X(H)$. We check using code of Rouse, Sutherland, and Zureick-Brown~\cite{RSZB_LAIOGFECOQ} that for each $r$ and each maximal subgroup $H \subset \GL_2(\bbZ/12\bbZ)$ there exists a prime $p$ as above such that $j(E_{12, r})$ is not the $j$-invariant of an $\bbF_p$-rational point on $X(H)$ (specifically we use the intrinsic \texttt{GL2jInvariantTest} in \cite[\texttt{gl2groups.m}]{RSZB_Electronic}). The claim follows from \thref{lemma:galoisrep}(ii).
		\end{proof}
	
	\subsection{Further points and curves on \texorpdfstring{$Z(12, r)$}{Z(12, r)}}
		Since the surfaces $Z(12, r)$ are of general type for $r \neq 1$ we would expect by the Bombieri--Lang conjecture that they contain at most finitely many curves of genus $0$ or $1$. In Tables~\ref{table:curves_Z125}--\ref{table:curves_Z1211} we record examples of curves of genus $0$ or $1$ on $Z(12, r)$ (where $r \neq 1$) which do not map to a Hecke correspondence on $Z(1)$ and do not lie above the loci where $j=0,1728,\infty$ or $j'= 0,1728,\infty$.

		In Tables~\ref{table:curves_Z125}--\ref{table:curves_Z1211} an irreducible curve $C \subset Z(12, r)$ is recorded by its image $C^+$ under the map $Z(12, r) \to W(12, r)$. The column ``genus'' records the genus of $C$. If $C$ is not geometrically irreducible we write $g^*$ where $g$ is the genus of a geometrically irreducible component. The column ``has $\bbQ$-point'' records whether a smooth projective curve birational over $\bbQ$ to $C$ has a rational point.

		\begin{table}[htbp]
			\centering 
			\begin{adjustbox}{width=\textwidth}
        \begin{tabular}{>{\centering\arraybackslash}p{5.5cm}|c|c|c|c|c}
          Equation for $C^+$                         & Genus of $C$	&{\makecell[c]{Primes where $C$ is \\locally insoluble}} & Has $\bbQ$-point	& {\makecell[c]{LMFDB label of \\the Jacobian of $C$}} & Rank \\
          \hline
          $u \pm v - 1$                              & $0$          & $\{3, \infty\} $                                       &No                &                                                      &     \\
          $v^2 + 1$                                  & $0^*$        &                                                        &No                &                                                      &     \\
          $-u^2 + v^2 + 2$                           & $0^*$        &                                                        &No                &                                                      &     \\
          $u^3 + 3u^2 - (v^2 + 10) u - 3 v^2 + 6$    & $1$          &                                                        &Yes               & \LMFDBLabel{240.a3}                                  & $1$ \\
          $-u + 2 \pm u  v \mp v + v^2$              & $1$          &                                                        &Yes               & \LMFDBLabel{32.a4}                                   & $0$ \\
          $u^2 - 2  u - v^2$                         & $1$          & $\{2, \infty\}$                                        &No                & \LMFDBLabel{32.a3}                                   & $0$ \\
          $u^2 + 4  u - 8 - v^2$                     & $1$          & $\{3, \infty\}$                                        &No                & \LMFDBLabel{45.a6}                                   & $0$ \\
          $u^2 - 2  u + 4 + 2  u  v - 2  v + v^2$    & $1^*$        &                                                        &No                &                                                      &     \\
          $u^2 - 2  u + 2 \pm 2  u  v \mp 2 v + v^2$ & $1^*$        &                                                        &No                &                                                      &     \\
        \end{tabular}
			\end{adjustbox}
      \caption{Examples of curves $C \subset Z(12, 5)$ of genus $0$ or $1$ without a known moduli interpretation.} \label{table:curves_Z125}
		\end{table}

		\begin{table}[htbp]
			\centering 
			\begin{adjustbox}{width=\textwidth}
        \begin{tabular}{>{\centering\arraybackslash}p{5.5cm}|c|c|c|c|c}
          Equation for $C^+$                  & Genus of $C$ &{\makecell[c]{Primes where $C$ is \\locally insoluble}} & Has $\bbQ$-point & {\makecell[c]{LMFDB label of \\the Jacobian of $C$}} & Rank \\
          \hline
          $u^2 + 3 u - 2 v^2 + 2$             & $0$          &                                                        & Yes              &                                                      &     \\
          $u + 2$                             & $0^*$        &                                                        & No               &                                                      &     \\
          $u^2 + 2u - 3 - v^2$                & $0^*$        &                                                        & No               &                                                      &     \\
          $u + 2 \pm 2v$                      & $1$          &                                                        & Yes              & \LMFDBLabel{2898.n2}                                 & $1$ \\
          $u + 2 \pm v$                       & $1$          &                                                        & Yes              & \LMFDBLabel{1224.h1}                                 & $1$ \\
          $u^2 - 2 u + 3 v^2 + 1$             & $1^*$        &                                                        & No               &                                                      &     \\
          $u^2 - 2 u + v^2 + 1$               & $1^*$        &                                                        & No               &                                                      &     \\
          $u^2 + 2 u + v^2 - 1$               & $1$          & $\{2, \infty\}$                                        & No               & \LMFDBLabel{64.a3}                                   & $0$ \\
          $u^2 - u + 2 v^2$                   & $1$          & $\{11, \infty\}$                                       & No               & \LMFDBLabel{7920.ba1}                                & $1$ \\
          $u^2 \pm u v + u + v^2 \pm 2 v + 1$	& $1^*$        &                                                        & No               &                                                      &     \\
        \end{tabular}
			\end{adjustbox}
      \caption{Examples of curves $C \subset Z(12, 7)$ of genus $0$ or $1$ without a known moduli interpretation.} \label{table:curves_Z127}
		\end{table}

		\begin{table}[htbp]
			\centering 
			\begin{adjustbox}{width=\textwidth}
        \begin{tabular}{>{\centering\arraybackslash}p{5.5cm}|c|c|c|c|c}
          Equation for $C^+$                            & Genus of $C$ &{\makecell[c]{Primes where $C$ is\\locally insoluble}} & Has $\bbQ$-point & {\makecell[c]{LMFDB label of \\the Jacobian of $C$}} & Rank \\
          \hline
          $(9 v^4 + 30 v^2 + 5) u - 45 v^4 - 6 v^2 + 7$	& $0$          &                                                        & Yes              &                                                      &     \\
          $u \pm 2v + 1$                                & $1$          &                                                        & Yes              & \LMFDBLabel{32.a4}                                   & $0$ \\
        \end{tabular}
			\end{adjustbox}
      \caption{Examples of curves $C \subset Z(12, 11)$ of genus $0$ or $1$ without a known moduli interpretation.} \label{table:curves_Z1211}
		\end{table}
		
		\begin{remark}
			It is natural to ask for which pairs $(N, r)$ there exist pairs of non-isogenous $(N, r)$-congruent elliptic curves $\mathcal{E}/\bbQ(t)$ and $\mathcal{E}'/\bbQ(t)$?
			
			In view of this question it would be interesting to determine whether or not there exists an example of a genus $0$ curve (which is birational over $\bbQ$ to $\bbP^1$) on $Z(12, 5)$ which does not lie above $j = \infty$, $j' = \infty$, or a Hecke correspondence on $Z(1)$. 
		\end{remark}

		\begin{remark}
			It is possible to find examples of curves of genus $0$ and $1$ on $Z(12, r)$ which do not arise as the image of a map $Y_0(m) \to Z(12, r)$ as described in Section~\ref{sec:X0m_on_Z} but which nevertheless map to the graph of a Hecke correspondence on $Z(1)$. 
			
			Examples are furnished by the genus $0$ curve $u + 2 = 0$ on $Z(12, 1)$ and the genus $1$ curve $u v^2 - v^2 - 2 = 0$ on $Z(12, 11)$. Points on these curves give rise to pairs of $12$-congruent elliptic curves which are $7$ and $5$-isogenous respectively.
			
			We do not attempt to classify all such examples. 
		\end{remark}

		When $r \neq 1$ the surfaces $Z(12, r)$ are of general type \cite{KS_MDQS} so we expect by the Bombieri--Lang conjecture that there are at most finitely many rational points on $Z(12, r)$ which do not lie on a known curve of genus $0$ or $1$. In total we were able to find 
		\begin{itemize}
		 \item $248$ points on $Z(12, 5)$ which correspond to $61$ pairs of $(12, 5)$-congruent elliptic curves, 
		 \item $232$ points on $Z(12, 7)$ which correspond to $58$ pairs of $(12, 7)$-congruent elliptic curves, and 
		 \item $8$ points on $Z(12, 11)$ which correspond to $2$ pairs of $(12, 11)$-congruent elliptic curves. 
		\end{itemize}
		These points are recorded (together with the corresponding congruent elliptic curves) in the electronic data attached to this paper~\cite{ME_ELECTRONIC_Z12}.

	\subsection{Examples of \texorpdfstring{$12$}{12}-congruent elliptic curves with small conductor} \label{sec:smallconductor}
				
		In Table \ref{table:small_ex} we give examples of non-isogenous $(12, r)$-congruent elliptic curves which have small conductor. For each $r \neq 1$ we list all examples which we found (up to simultaneous quadratic twist) with conductor $\leq 10^{10}$. These examples were constructed by searching for rational points on $Z(12, r)$. When $r = 1$ we were able to find too many examples to record here, the $18$ examples we found where both $E$ and $E'$ are contained in the LMFDB are recorded in the electronic data~\cite{ME_ELECTRONIC_Z12}.

		We do not include examples where $E$ and $E'$ are quadratic twists since many such examples may be constructed for $r = 7, 11$ by searching for points on the modular curves in \cite[Lemmas~6.2 and 6.3]{F_COECAFNSMNGR}.

		The pair of $(12, 5)$-congruent elliptic curves \LMFDBLabel{60450.cx2} and \LMFDBLabel{60450.cw2} appearing in Table~\ref{table:small_ex} are notable since they give rise to $8$ points on $Z(12, 5)$. This corresponds to the fact that \LMFDBLabel{60450.cx2} admits a non-trivial $G_\bbQ$-equivariant automorphism of its $12$-torsion subgroup with power $1$.

    \begin{table}[H]
			\centering 
			\begin{adjustbox}{height=4.5cm}
			\begin{tabular}{c >{\centering\arraybackslash}p{2cm} >{\centering\arraybackslash}p{2cm} | c >{\centering\arraybackslash}p{2cm} >{\centering\arraybackslash}p{2cm}}	
				$r$  & $E$                            & $E'$                            & $r$  & $E$                            & $E'$                           \\
				\hline
				$7$  & \LMFDBLabel{735.d2}            & \LMFDBLabel{24255.bh1}          & $5$  & \texttt{16039350*}             & \texttt{176432850*}            \\ 
				$7$  & \LMFDBLabel{11830.x2}          & \LMFDBLabel{11830.m1}           & $7$  & \texttt{25860114*}             & \texttt{336181482*}            \\ 
				$5$  & \LMFDBLabel{60450.cx2}         & \LMFDBLabel{60450.cw2}          & $7$  & \texttt{28594930*}             & \texttt{2030240030*}           \\ 
				$5$  & \LMFDBLabel{76614.o1}          & \texttt{1302438*}               & $7$  & \texttt{29894410*}             & \texttt{29894410*}             \\ 
				$5$  & \LMFDBLabel{92510.w1}          & \LMFDBLabel{92510.i1}           & $5$  & \texttt{39160576*}             & \texttt{39160576*}             \\ 
				$5$  & \texttt{769515*}               & \texttt{45401385*}              & $7$  & \texttt{123160716*}            & \texttt{5788553652*}           \\ 
				$7$  & \texttt{838240*}               & \texttt{838240*}                & $7$  & \texttt{124010865*}            & \texttt{124010865*}            \\ 
				$5$  & \texttt{1991985*}              & \texttt{1991985*}               & $7$  & \texttt{159365310*}            & \texttt{159365310*}            \\ 
				$5$  & \texttt{4353790*}              & \texttt{4353790*}               & $11$ & \texttt{388152270*}            & \texttt{388152270*}            \\ 
				$11$ & \texttt{4976690*}              & \texttt{4976690*}               & $5$  & \texttt{1885555650*}           & \texttt{1885555650*}           \\ 
				$5$  & \texttt{5144914*}              & \texttt{5144914*}               & $5$  & \texttt{2016538368*}           & \texttt{2016538368*}           \\ 
				$5$  & \texttt{5355525*}              & \texttt{91043925*}              & $5$  & \texttt{2992233426*}           & \texttt{2992233426*}           \\ 
				$7$  & \texttt{10691890*}             & \texttt{10691890*}              & $7$  & \texttt{8900609640*}           & \texttt{8900609640*}           \\ 
				$5$  & \texttt{13286570*}             & \texttt{13286570*}              &      &                                &                                \\
			\end{tabular}
    \end{adjustbox}
    \caption{Known examples of pairs of geometrically non-isogenous $(12, r)$-congruent elliptic curves $E/\bbQ$ and $E'/\bbQ$ (up to simultaneous quadratic twist) of conductors $\leq 10^{10}$ where $r \neq 1$. We write $N*$ for an elliptic curve of conductor $N$ which does not appear in the LMFDB~\cite{lmfdb} (i.e., of conductor $N \geq 500\,000$). } \label{table:small_ex}
		\end{table}

    \providecommand{\bysame}{\leavevmode\hbox to3em{\hrulefill}\thinspace}
    \providecommand{\MR}{\relax\ifhmode\unskip\space\fi MR }
    \providecommand{\MRhref}[2]{%
      \href{http://www.ams.org/mathscinet-getitem?mr=#1}{#2}
    }
    \providecommand{\href}[2]{#2}

\end{document}